\documentclass[a4paper,12pt]{article}
\usepackage[english]{babel}
\usepackage[latin1]{inputenc}
\usepackage[T1]{fontenc}
\usepackage{graphics}
\usepackage{booktabs}
\usepackage{amsmath}
\usepackage{inputenc}
\usepackage{amsthm}
\usepackage{amssymb}
\usepackage{setspace}
\usepackage{array}
\usepackage{graphicx}
\usepackage{subfig}
\usepackage{varioref}
\usepackage[blocks]{authblk}
\usepackage[font=footnotesize,format=hang]{caption}
\usepackage[a4paper,top=3cm,bottom=3cm,left=2.5cm,right=2.5cm,%
bindingoffset=5mm]{geometry}

\theoremstyle{plain}
\newtheorem{proposition}{Proposition}

\theoremstyle{remark}
\newtheorem{observation}{Observation}

\title{\textbf{A new magnitude-dependent ETAS model for earthquakes}}
\date{}
\usepackage{authblk}

\author[1]{Ilaria Spassiani}
\author[1]{Giovanni Sebastiani}
\affil[1]{Department of Mathematics ``Guido Castelnuovo'', Sapienza University of Rome, Rome, ITALY.}
\affil[2]{Istituto per le Applicazioni del Calcolo ``M. Picone'',
Consiglio Nazionale delle Ricerche, Rome, ITALY.}

\begin{document}
\maketitle

\begin{abstract}
\noindent We propose a new version of the ETAS model, which we also analyze theoretically. As for the standard ETAS model, we assume the Gutenberg-Richter law as a probability density function for background events' magnitude. Instead, the magnitude of triggered shocks is assumed to be probabilistically dependent on the triggering events' magnitude. To this aim, we propose a suitable probability density function. This function is such that, averaging on all triggering events' magnitudes, we obtain again the Gutenberg-Richter law. This ensures the validity of this law at any event's generation when ignoring past seismicity. The probabilistic dependence between the magnitude of triggered events' and the one of the corresponding triggering shock is motivated by some results of a statistical analysis of some Italian catalogues. The proposed model has been also theoretically analyzed here. In particular, we focus on the inter-event time which plays a very important role in the assessment of seismic hazard~\cite{molchan:primo}. Using the tool of the probability generating function and the Palm theory, we derive the density of inter-event time for small values.
\end{abstract}

\section{\textbf{Introduction}}
The Epidemic Type Aftershock Sequences (ETAS) is a model based on a specific branching process and represents a benchmark in statistical and mathematical seismology. The process generates the events along consecutive generations. At any step, each event of a given generation produces, independently of the others, its own cluster of events. The number of elements in a cluster is a Poisson random variable dependent on the generating event's magnitude through the \emph{productivity law}. Occurring times and magnitudes of each event in a cluster are independent of each other and are distributed according to the \emph{Omori-Utsu} and the \emph{Gutenberg-Richter laws}, respectively (for the above-mentioned three laws see equations~\eqref{eqn:GR},~\eqref{eqn:omori} and~\eqref{eqn:produttivita} in Appendix~\ref{sec:appc}). These times and magnitudes are also independent of the characteristics of past events. Similarly, given an event, the elements relative to its cluster are spatially distributed according to a particular model centered at it (see~\cite{ogata:secondo}). Here we will consider the process only regarding to the time, ignoring the spatial location of events. In particular, we will focus on the random variables relative to inter-event times, that are the times between consecutive shocks.

In the ETAS model, seismic events are divided into spontaneous (or background) and triggered ones. The characteristics of the branching process involved and the homogeneity of the background component imply that the total process is homogeneous, too. Therefore, we can deduce that inter-event times have all the same distribution. We will then consider the single inter-event time between consecutive shocks. The characteristics of this random variable play a very important role in the assessment of seismic hazard~\cite{molchan:primo}. In what follows, we will assume that $\lambda$ is the rate of the whole process of ''observable'' events and $\omega$ is the one corresponding to the ''observable'' events of the spontaneous component.
\begin{observation}
Given the events of any seismic catalogue, each of them is considered observable only if it has a magnitude larger than a certain cutoff value $m$. This threshold is known as the \emph{completeness magnitude}, that is the smallest value such that catalog events' magnitudes exceeding that value are distributed according to the Gutenberg-Richter model. The deviations from that model, when considering catalogue magnitudes smaller than $m$, can be interpreted in terms of wrong or missing events measurements.
\end{observation}

A very useful tool to study the distribution of the inter-event time is the probability generating function (PGF) of the random number $N(\tau)$ of observable events occurring in the time interval $[0,\tau]$, that is the quantity
\begin{equation}
\label{eqn:40}
G_{N(\tau)}(z)=\mathbb{E}[z^{N(\tau)}].
\end{equation}
The value of this function in zero gives the probability $\mathbb{P}(\tau)$ to have zero events in $[0,\tau]$. The inter-event time density is then obtained through the Palm equation by a double time derivative of $\mathbb{P}(\tau)$ and scaling. This has been done for the standard ETAS model by Saichev and Sornette~\cite{sornettesaichev:primo}. However, the standard ETAS model assumes that the magnitude of each event in any cluster is distributed according to the Gutenberg-Richter law, independently of any past event. Instead, it is assumed that the number of triggered events in a cluster depends on the triggering event's magnitude. We find intuitive and more realistic to also assume that the distribution of triggered events' magnitude in a cluster depends on the triggering events' magnitude. This is supported by evidences of an empirical study (see~\cite{mio:primo}). Therefore we propose a new ETAS version which assumes the above-cited hypothesis. Following the same technical tool of Saichev and Sornette, that is the above-mentioned PGF, we derive an expression for the inter-event time density for small values. In the following section we are going to compute the above-mentioned probability generating function corresponding to the proposed model.

\section{\textbf{Probability generating function for earthquake sequences}}
As previously specified, we are interested in $N(\tau)$, that is the number of observable events in $[0,\tau]$, which are of two types: spontaneous and triggered. Let's consider at first the triggered events. We can associate to each of these events the spontaneous shock from which the branching process reaches that event. Let's consider then all spontaneous events. The triggered events can be partitioned according to their associated spontaneous triggering shocks. The total number $N_1(\tau)$ of observable events triggered in $[0,\tau]$ will be the sum of observable events triggered by each spontaneous event. Let's consider now a partition of the infinite interval of times smaller than $\tau$ into small subintervals $I_k=[t_k,t_k+\Delta t]$ grouped in two categories. To the first (second) category belong subintervals of kind $I_{k_1}\subset(-\infty,0)$ ($I_{k_2}\subset[0,\tau]$). The probability generating function of the number $N_1(\tau)$ can be obtained as the limit, with respect to $\Delta t$, of the corresponding probability generating function relative to the subinterval partition.

Recall now both the property of the PGF expressed in equation~\eqref{prop:2}, Appendix~\ref{sec:appb}, and the fact that there's independence between the numbers of events, triggered by the spontaneous ones belonging to disjoint subintervals. We can therefore deduce that the PGF of $N_1(\tau)$ is the limit of the product between the respective probability generating functions relative to the different subintervals. Let's consider the PGF of the number of events in $[0,\tau]$, triggered by spontaneous shocks contained in a generic subinterval. Thanks to the properties of the probability generating function, this PGF is obtained by composing the PGF of the number of background events occurred in that subinterval, with the PGF relative to the number of events in $[0,\tau]$ triggered by just one of the considered spontaneous events (see Proposition~\eqref{prop:1}, Appendix~\ref{sec:appb}).

Furthermore, in order to obtain the probability generating function of the total number $N(\tau)$ of observable events, we have to take into account spontaneous events, too. Obviously, the background events to be considered are only those which fall in subintervals of the second type. In this case, for each of these events we will count all the shocks triggered by it, plus the event itself. This last is equivalent to multiply by $z$ the probability generating function of the number of events triggered in $[0,\tau]$ by a single spontaneous event of magnitude $m'$, itself belonging to $[0,\tau]$. More precisely, since this event has to be observable, the multiplication by $z$ will be done only if its magnitude is larger than the threshold $m$. Therefore, instead of the multiplication by $z$, we will multiply by
\begin{equation}
\label{eqn:perz}
f(z,m',m)=U(m'<m)+zU(m'\ge m),
\end{equation}
where $U(x>0)$ is the indicator function having value one where $x>0$, zero otherwise.

Let's now compute explicitly the probability generating function $\Omega(\cdot;\tau)$ of $N(\tau)$. Let's consider the number of events with magnitude larger than $m$, triggered in a generic interval $[0,\nu]$ by a background event with magnitude $m'$ occurred in a time $-s$, where $s\ge0$. In what follows, we will indicate with $G(\cdot;s,\nu,m')$ the PGF of such a number. Recalling that the random variable associated with the number of spontaneous events in a generic subinterval of length $\Delta t$ follows a Poisson distribution with parameter $\omega\Delta t$, for the previous reasonings we have
\begin{align}
\label{eqn:omegagrande}
\Omega(z;\tau)=&\lim_{\Delta t\to0}\prod_{I_{k_1}}\exp\Biggl\{-\omega\Delta t[1-G(z;s,\tau,m')]\Biggr\}\notag\\
&\cdot\prod_{I_{k_2}}\exp\Biggl\{-\omega\Delta t[1-f(z,m',m)G(z;0,\tau-w,m')]\Biggr\}\notag\\
=&\exp \Biggl\{-\omega\biggl[\biggl(\int_{-\infty}^0[1-G(z;-t^\prime,\tau,m^\prime)]dt^\prime\notag\\
&+\int_0^{\tau}[1-f(z,m',m)G(z;0,\tau-t^\prime,m^\prime)]dt^\prime\biggr)   \biggr]\Biggr\}\notag\\
=&\exp \Biggl\{-\omega\biggl[\biggl(\int_{0}^\infty[1-G(z;t,\tau,m^\prime)]dt\notag\\
&+\int_0^{\tau}[1-f(z,m',m)G(z;0,t,m^\prime)]dt\biggr)   \biggr]      \Biggr\},
\end{align}
where $-s\in I_{k_1}$ and $w\in I_{k_2}$.

The formulae obtained till here are relative to observable events with magnitude larger than $m$, triggered by a spontaneous shock with a fixed but generic magnitude equal to $m'$. Nevertheless, we have to introduce another threshold value $m_0$, called \emph{reference magnitude}, which corresponds to the minimum value for an event to be able to trigger their own offsprings~\cite{sornhelmn:primo}. We therefore have to substitute in~\eqref{eqn:omegagrande} the generating function $G(z;\cdot,\cdot,m')$ with
\[
\int_{m_0}^\infty G(z;\cdot,\cdot,m')p(m')dm',
\]
where $p(m')$ is the Gutenberg-Richter law (equation~\eqref{eqn:GR} in Appendix~\ref{sec:appc}).
\begin{observation}
Reference magnitude is usually set less than or equal to the completeness one~\cite{sornettewerner:primo}. Indeed, even if the opposite inequality is verified, we could however reduce to the previous case by excluding the data relative to observable events with magnitude less than $m_0$.
\end{observation}

The function $p(\cdot)$ is universally adopted to model the probability density function of a generic event's magnitude, if we don't take into account the characteristics of its previous seismicity. In the specific case we are considering background events; we assume that their magnitude is also modeled by the function $p(\cdot)$. Obviously it holds
\[
\int_{m_0}^\infty p(m')dm'=1.
\]
Setting
\begin{align}
\label{eqn:b12}
L(z;\tau)=&\int_{m_0}^\infty dm'p(m')\biggl(\int_0^\infty[1-G(z;t,\tau,m')]dt\notag\\
             &+\int_0^\tau[1-f(z,m',m)G(z;0,t,m')]dt\biggr),
\end{align}
we can write, with abuse of notation,
\begin{equation}
\label{eqn:b11}
\Omega(z;\tau)=e^{-\omega L(z;\tau)}.
\end{equation}
Let's substitute now the above-cited function $f(z,m',m)$ (equation~\eqref{eqn:perz}) in the last expression of~\eqref{eqn:b12}. We obtain
\begin{align}
\label{eqn:b18}
L(z;\tau)=&\int_{m_0}^\infty dm'p(m')\biggl(\int_{0}^\infty[1-G(z;t,\tau,m')]dt\notag\\
           &+\int_{0}^\tau\Bigl[1-[U(m'<m)+zU(m'\ge m)]G(z;0,t,m')\Bigr]dt\biggr)\notag\\
           =&\int_0^\infty\biggl[1-\int_{m_0}^\infty G(z;t,\tau,m')p(m')dm'\biggr]dt\notag\\
           &+\int_{0}^\tau\biggl[1-\int_{m_0}^m G(z;0,t,m')p(m')dm'\notag\\
           &-z\int_m^\infty G(z;0,t,m')p(m')dm'\biggr]dt\notag\\
           =&\int_0^\infty\biggl[1-\int_{m_0}^\infty G(z;t,\tau,m')p(m')dm'\biggr]dt\notag\\
           &+\int_0^\tau\biggl[1-\int_{m_0}^\infty G(z;0,t,m')p(m')dm'\notag\\
           &+(1-z)\int_m^\infty G(z;0,t,m')p(m')dm'\biggr]dt.
\end{align}

At this point, we want to compute explicitly the probability generating function $\Omega(\cdot;\tau)$. To do that, we find an expression for $G(z;t,\tau,m')$ that is the PGF of the number of observable events in $[0,\tau]$ triggered by a fixed spontaneous shock with magnitude $m'$ occurred in $-t$, where $t\ge0$. Let's consider now any observable event triggered in $[0,\tau]$ belonging to any generation from the second one on and let $N_2(\tau)$ denote their total number. As before, we can associate it with a first generation triggered shock occurred between the spontaneous one and the event itself. The number $N_2(\tau)$ will be the sum of observable events triggered by each first generation triggered shock. Thanks to the auto-similarity of the process, the above-mentioned probability generating function $G(z;\cdot,\cdot,\cdot)$ is the same for both a background event and a triggered one. Therefore, as previously done, we obtain
\begin{align}
\label{eqn:altra}
G(z;t,\tau,&m')=\lim_{\Delta t\to0}\prod_{I_{k_1}\subset[-t,0]}\exp\Biggl\{\Phi(-s+t)\Delta t\rho(m')[G(z;s,\tau,m'')-1]\Biggr\}\notag\\
                  &\cdot\prod_{I_{k_2}}\exp\Biggl\{\Phi(w+t)\Delta t\rho(m')[f(z,m'',m)G(z;0,\tau-w,m'')-1]\Biggr\},
\end{align}
where $-s\in I_{k_1}$, $w\in I_{k_2}$ and $m''$ is the magnitude of the first generation triggered shock. Recall that $\Phi(\cdot)$ is the Omori-Utsu law and $\rho(\cdot)$ is the productivity one (respectively equations~\eqref{eqn:omori} and~\eqref{eqn:produttivita} in Appendix~\ref{sec:appc}). Again, we will have to average with respect to the first generation triggered shock's magnitude. To this end, we assume that the latter magnitude, here indicated with $m''$ conditioned on the background event's magnitude $m'$, follows the model $\mathbb{P}(m''|m')=p(m''|m')$, where the proposal for the last probability density  is given in Appendix~\ref{sec:appa}. Obviously it must hold
\[
\int_{m_0}^\infty p(m''|m')dm''=1.
\]

Any first generation triggered event considered is also a triggering one. Then, its magnitude will have to be integrated from $m_0$ to infinity. Recalling equation~\eqref{eqn:perz}, from~\eqref{eqn:altra} we obtain, with abuse of notation, the following expression:
\begin{align}
\label{eqn:b13}
G(z;t,\tau,m')=&\exp\Biggl\{-\rho(m')\int_{m_0}^\infty dm''p(m''|m')\biggl[\int_{-t}^0\Phi(t+x)[1\notag\\
                  &-G(z;-x,\tau,m'')]dx+\int_0^\tau\Phi(t+x)[1\notag\\
                  &-f(z,m'',m)G(z;0,\tau-x,m'')]dx\biggr]\Biggr\}\notag\\
                  =&\exp\Biggl\{-\rho(m')\int_{m_0}^\infty dm''p(m''|m')\biggl[\int_{-t}^0\Phi(t+x)dx\notag\\
                  &-\int_{-t}^0\Phi(t+x)G(z;-x,\tau,m'')dx+\int_0^\tau\Phi(t+x)dx\notag\\
                  &-f(z,m'',m)\int_0^\tau \Phi(t+x)G(z;0,\tau-x,m'')dx\biggr]\Biggr\}\notag\\
                  =&\exp\Biggl\{-\rho(m')\int_{m_0}^\infty dm''p(m''|m')\biggl[\int_{0}^{t}\Phi(x)dx\notag\\
                  &-\int_{0}^{t}\Phi(t-x)G(z;x,\tau,m'')dx+\int_t^{t+\tau}\Phi(x)dx\notag\\
                  &-f(z,m'',m)\int_0^{\tau} \Phi(t+\tau-x)G(z;0,x,m'')dx\biggr]\Biggr\}\notag\\
                  =&\exp\Biggl\{-\rho(m')\int_{m_0}^\infty dm''p(m''|m')[b(t+\tau)\notag\\
                  &-(\Phi(\cdot) \otimes G(z;\cdot,\tau,m''))(t)\notag\\
                  &-f(z,m'',m)(\Phi(t+\cdot)\otimes G(z;0,\cdot,m''))(\tau)]\Biggr\}\notag\\
                  =&\exp\Biggl\{-\rho(m')\biggl[b(t+\tau)\notag\\
                  &-\int_{m_0}^\infty(\Phi(\cdot) \otimes G(z;\cdot,\tau,m''))(t)p(m''|m')dm''\notag\\
                  &-\int_{m_0}^\infty(\Phi(t+\cdot)\otimes G(z;0,\cdot,m''))(\tau)p(m''|m')dm''\notag\\
                  &+(1-z)\int_{m}^\infty(\Phi(t+\cdot)\otimes G(z;0,\cdot,m''))(\tau)p(m''|m')dm''\biggr]\Biggr\},\notag\\
\end{align}
where $\otimes$ is the convolution operator and
\begin{equation}
\label{eqn:b}
b(t)=\int_0^{t}\Phi(x)dx.
\end{equation}

At this point, in order to find the probability of having zero events in $[0,\tau]$, we will have to substitute equation~\eqref{eqn:b13} in~\eqref{eqn:b18} opportunely, reminding that
\[
\mathbb{P}(\tau)=\Omega(0;\tau)=e^{-\omega L(0;\tau)}.
\]
In what follows we will work separately on two cases. First of all on the particular case $m=m_0$, then on the general one $m\ge m_0$.

\section{\textbf{The case $m=m_0$}}
In this section, we want to focus on the particular case in which $m=m_0$. We will focus on this case because we can find a closed-form expression for the probability $\mathbb{P}(\tau)$ of having zero events in $[0,\tau]$, for small values of $\tau$. Than, by using the Palm theory, we derive a closed-form expression for the density of inter-event time. The relevance of the case $m=m_0$ can be understood by the following reasoning. On one side, the reference magnitude $m_0$ is connected to the intrinsic physical features of soil and remains almost constant in a sufficiently long time frame. On the other side, the completeness magnitude is related to the number of seismic stations and their sensitivity. Then, the latter magnitude varies with seismic stations density.
Hoping for a gradual improvement of technology both in instrumental sensitivity and in the number of seismic stations available, we can certainly suppose to reduce $m$ till it will be just equal to $m_0$.

Setting $m=m_0$ in the last expression of formula~\eqref{eqn:b18} and recalling~\eqref{eqn:b11}, we find that the probability generating function $\Omega(z;\tau)$ relative to the total number of observable events in $[0,\tau]$ is
\begin{align}
\label{eqn:sottocaso1}
\Omega(z;\tau)=&\exp\Biggl\{-\omega\biggl[\int_0^\infty\Bigl[1-\int_{m_0}^\infty G(z;t,\tau,m')p(m')dm'\Bigr]dt\notag\\
&+\int_0^\tau\Bigl[1-z\int_{m_0}^\infty G(z;0,t,m')p(m')dm'\Bigr]dt\biggr]\Biggr\}.
\end{align}
In this particular case the function $L(0;\tau)$ becomes
\begin{align}
\label{eqn:sottocaso4}
L(0;\tau)&=\int_0^\infty\Bigl[1-\int_{m_0}^\infty G(0;t,\tau,m')p(m')dm'\Bigr]dt+\int_0^\tau dt\notag\\
&=\int_0^\infty\Bigl[1-\int_{m_0}^\infty G(0;t,\tau,m')p(m')dm'\Bigr]dt+\tau\notag\\
&=\int_0^\infty\overline{N}_{m=m_0}(t,\tau)dt+\tau,
\end{align}
with
\begin{equation}
\label{eqn:sottocaso10}
\overline{N}_{m=m_0}(t,\tau)=1-\int_{m_0}^\infty G(0;t,\tau,m')p(m')dm'.
\end{equation}

Let's compute then the probability generating function $G(0;t,\tau,m')$. Setting $z=0$ and $m=m_0$ in the last expression of~\eqref{eqn:b13} and considering the approximation of the exponential function, we obtain
\begin{align}
\label{eqn:sottocaso2}
G(0;t,\tau,m')=&\exp\Biggl\{-\rho(m')\biggl[b(t+\tau)\notag\\
&-\int_{m_0}^\infty(\Phi(\cdot) \otimes G(0;\cdot,\tau,m''))(t)p(m''|m')dm''\biggr]\Biggr\}\notag\\
\approx&1-\rho(m')\biggl[b(t+\tau)\notag\\
&-\int_{m_0}^\infty(\Phi(\cdot) \otimes G(0;\cdot,\tau,m''))(t)p(m''|m')dm''\biggr].
\end{align}
This approximation in Taylor series around the point zero makes sense since the exponential argument is zero when $\tau=0$. In fact, $G(0;\cdot,0,m'')=1$.
If we substitute the previous formula in~\eqref{eqn:sottocaso10}, we get
\begin{align}
\label{eqn:sottocaso7}
\overline{N}_{m=m_0}(t,\tau)\approx&1-\int_{m_0}^\infty dm'p(m')\biggl\{1-\rho(m')\Bigl[b(t+\tau)\notag\\
&-\int_{m_0}^\infty(\Phi(\cdot) \otimes G(0;\cdot,\tau,m''))(t)p(m''|m')dm''\Bigr]\biggr\}\notag\\
=&1-1+b(t+\tau)\int_{m_0}^\infty p(m')\rho(m')dm'\notag\\
&-\int_{m_0}^\infty dm'p(m')\rho(m')\int_{m_0}^\infty(\Phi(\cdot)\otimes G(0;\cdot,\tau,m''))(t)p(m''|m')dm''\notag\\
=&nb(t+\tau)-\int_{m_0}^\infty dm'p(m')\rho(m')\int_{m_0}^\infty dm''p(m''|m')\notag\\
&\cdot\int_0^t\Phi(x)G(0;t-x,\tau,m'')dx\notag\\
&\text{switch the order of integration}\notag\\
=&nb(t+\tau)-\int_0^tdx\,\Phi(x)\int_{m_0}^\infty dm''G(0;t-x,\tau,m'')\notag\\
&\cdot\int_{m_0}^\infty p(m')\rho(m')p(m''|m')dm'\notag\\
=&nb(t+\tau)-\int_0^tdx\,\Phi(x)\int_{m_0}^\infty G(0;t-x,\tau,m'')I_A(m'')dm'',
\end{align}
with
\begin{equation}
\label{eqn:Agrande}
I_A(m'')=\int_{m_0}^\infty p(m')\rho(m')p(m''|m')dm'
\end{equation}
and
\begin{align}
\label{eqn:enne}
n&=\int_{m_0}^\infty p(m')\rho(m')dm'\notag\\
&=\int_{m_0}^\infty\beta e^{-\beta(m'-m_0)}\kappa e^{a(m'-m_0)}dm'\notag\\
&=\beta\kappa\int_{m_0}^\infty e^{-(\beta-a)(m'-m_0)}dm'\notag\\
&=\frac{\beta\kappa}{-(\beta-a)} e^{-(\beta-a)(m'-m_0)}\bigg|^\infty_{m_0}\notag\\
&=\frac{\beta\kappa}{\beta-a}.
\end{align}
Equation~\eqref{eqn:enne} has been obtained using the formulae~\eqref{eqn:GR} and~\eqref{eqn:produttivita} in Appendix~\ref{sec:appc}.

Let's assume now that the following condition is verified:
\begin{align*}
\frac{\int_{m_0}^\infty p(m')\rho(m')p(m''|m')dm'}{\int_{m_0}^\infty p(m')\rho(m')dm'}&=p(m'')\\
&\Updownarrow\\
I_A(m'')&=np(m'').
\end{align*}
\begin{observation}
\label{oss:invperpcond}
This condition represents the invariance of the Gutenberg-Richter law, weighted by $\rho(m')$, with respect to the conditional probability density function assumed valid for triggered events' magnitude. It is a focal point of our study. Although we assume the existence of a density of triggered events' magnitude different from the Gutenberg-Richter law, at the same time the above condition justifies its validity. More precisely, the condition we are imposing tells us that, if we average the conditional density $p(m''|m')$ over all the possible triggering events' magnitudes $m'$, distributed according to the Gutenberg-Richter law and taking into account the productivity model, we obtain again the above-cited law. This shows that first generation events' magnitude follows the Gutenberg-Richter law as assumed for background shocks. In addition, this property will be obviously true for events' magnitude of any generation. Hence, iterating the reasoning, a generic event's magnitude, whatever generation it belongs to, will be distributed according to the Gutenberg-Richter law.
\end{observation}
Under this assumption, equation~\eqref{eqn:sottocaso7} becomes
\begin{align}
\label{eqn:sottocaso8}
\overline{N}_{m=m_0}(t,\tau)\approx& nb(t+\tau)-n\int_0^tdx\,\Phi(x)\int_{m_0}^\infty G(0;t-x,\tau,m'')p(m'')dm''\notag\\
=&n\bigl[b(t+\tau)-b(t)\bigr]+n\biggl[b(t)\notag\\
&-\int_0^tdx\,\Phi(x)\int_{m_0}^\infty G(0;t-x,\tau,m'')p(m'')dm''\biggr]\notag\\
=&n\bigl[b(t+\tau)-b(t)\bigr]+n\int_0^tdx\,\Phi(x)\biggl[1\notag\\
&-\int_{m_0}^\infty G(0;t-x,\tau,m'')p(m'')dm''\biggr]\notag\\
=&n\bigl[b(t+\tau)-b(t)\bigr]+n\int_0^t\Phi(x)\overline{N}_{m=m_0}(t-x,\tau)dx,
\end{align}
where in the last equality we've used formula~\eqref{eqn:sottocaso10}.
For simplicity of notation, in what follows we will substitute approximation symbol with equal sign.

In order to obtain the expression to substitute in~\eqref{eqn:sottocaso4}, let's compute the temporal integral of both the first and the last members of equation~\eqref{eqn:sottocaso8}. We get
\begin{align*}
\int_0^\infty\overline{N}_{m=m_0}(t,\tau)dt=&n\int_{0}^\infty\bigl[b(t+\tau)-b(t)\bigr]dt +n\int_0^\infty\int_0^t\Phi(t-x)\overline{N}_{m=m_0}(x,\tau)dxdt\\
&\text{switch the last two temporal integrals}\\
=&n\int_{0}^\infty\bigl[b(t+\tau)-b(t)\bigr]dt +n\int_0^\infty\int_x^\infty\Phi(t-x)\overline{N}_{m=m_0}(x,\tau)dtdx\\
=&n\int_{0}^\infty\bigl[b(t+\tau)-b(t)\bigr]dt +n\int_0^\infty\overline{N}_{m=m_0}(x,\tau)dx,
\end{align*}
where the last equality is obtained because, recalling~\eqref{eqn:omori}, it obviously holds
\begin{equation}
\label{eqn:omdens}
\int_x^\infty\Phi(t-x)dt=\int_0^\infty\Phi(y)dy=1.
\end{equation}
From~\eqref{eqn:omdens} we then have
\begin{align}
\label{eqn:altraa}
\int_0^\infty\overline{N}_{m=m_0}(t,\tau)dt&=\frac{n}{1-n}\int_0^\infty\bigl[b(t+\tau)-b(t)\bigr]dt\notag\\
&=\frac{n}{1-n}\int_0^\infty\biggl[1-\frac{c^\theta}{(c+t+\tau)^\theta}-1+\frac{c^\theta}{(c+t)^\theta}\biggr]dt\notag\\
&=\frac{nc^\theta}{1-n}\int_0^\infty\biggl[\frac{1}{(c+t)^\theta}-\frac{1}{(c+t+\tau)^\theta}\biggr]dt\notag\\
&=\frac{nc^\theta}{(1-n)(1-\theta)} \Bigl[(c+t)^{1-\theta}-(c+t+\tau)^{1-\theta}\Bigr]\bigg|_0^\infty\notag\\
&=\frac{nc^\theta}{(1-n)(1-\theta)} \Bigl[(c+\tau)^{1-\theta}-c^{1-\theta}\Bigr].
\end{align}
The second equality follows from the explicit form of the function $b(t)$ defined in~\eqref{eqn:b}:
\begin{align*}
b(t)&=\int_0^t\frac{\theta c^\theta}{(c+x)^{1+\theta}}dx\\
&=\theta c^\theta\int_0^t(c+x)^{-1-\theta}\\
&=-c^\theta(c+x)^{-\theta}\big|_0^t\\
&=1-\frac{c^\theta}{(c+t)^\theta}.
\end{align*}
On the other hand, the last equality in~\eqref{eqn:altraa} is obtained since the limit of the function $(c+t)^{1-\theta}-(c+t+\tau)^{1-\theta}$ is, for $t$ tending to infinity, always zero. In fact, for $\theta>1$ it is obvious; if instead $1>\theta>0$ we have
\begin{align*}
\lim_{t\to\infty}\bigl[(c+t)^{1-\theta}-(c+t+\tau)^{1-\theta}\bigr]=&\lim_{t\to\infty}(c+t)^{1-\theta}\biggl[1-\Bigl(1+\frac{\tau}{c+t}\Bigr)^{1-\theta}\biggr]\\
=&\lim_{t\to\infty}\frac{1-\Bigl(1+\frac{\tau}{c+t}\Bigr)^{1-\theta}}{(c+t)^{\theta-1}}\\
&\text{apply De L'H\^{o}pital's rule}\\
=&\lim_{t\to\infty}\frac{(1-\theta)\Bigl(1+\frac{\tau}{c+t}\Bigr)^{-\theta}\frac{\tau}{(c+t)^2}}{(\theta-1)(c+t)^{\theta-2}}\\
=&\lim_{t\to\infty}-\frac{\Bigl(1+\frac{\tau}{c+t}\Bigr)^{-\theta}\tau}{(c+t)^\theta}\\
=&0.
\end{align*}
At this point we can substitute equation~\eqref{eqn:altraa} in~\eqref{eqn:sottocaso4}:
\[
L(0;\tau)=\frac{n}{(1-n)(1-\theta)} \Bigl[c^\theta(c+\tau)^{1-\theta}-c\Bigr]+\tau.
\]
Introducing the functions
\begin{equation}
\label{eqn:funa}
a(t)=1-b(t)=\int_t^\infty\Phi(x)dx=\frac{c^{\theta}}{(c+t)^{\theta}}
\end{equation}
and
\begin{equation}
\label{eqn:funagrande}
A(\tau)=\int_0^\tau a(x)dx=\frac{c^{\theta}}{1-\theta}(c+t)^{1-\theta}\big|_0^{\tau}=\frac{\bigl[c^{\theta}(\tau+c)^{1-\theta}-c\bigr]}{1-\theta},
\end{equation}
we get
\begin{align}
\label{eqn:sottocaso9}
\mathbb{P}(\tau)&=e^{-\omega L(0;\tau)}\notag\\
&=\exp\Biggl\{-\omega\tau-\frac{n\omega}{1-n}A(\tau)\Biggr\}.
\end{align}
As already specified, using the Palm theory we obtain the density $F_{\tau}(\tau)$ relative to the inter-event time in the following way (see~\cite{daleyverejones:primo}):
\begin{equation}
\label{eqn:intensintert}
F_{\tau}(\tau)=\frac{1}{\lambda}\frac{d^2\mathbb{P}(\tau)}{d\tau^2}.
\end{equation}
Let's compute than the second derivative of the above-mentioned probability. We have
\[
\frac{d\mathbb{P}(\tau)}{d\tau}=\mathbb{P}(\tau)\Bigl[-\omega-\frac{n\omega}{1-n}a(\tau)\Bigr],
\]
in fact
\begin{equation}
\label{eqn:derivA}
\frac{d}{d\tau}A(\tau)=\frac{d}{d\tau}\int_0^{\tau}a(x)dx=a(\tau).
\end{equation}
Hence,
\begin{align*}
\frac{d^2\mathbb{P}(\tau)}{d\tau^2}=&\mathbb{P}(\tau)\Bigl[-\omega-\frac{n\omega}{1-n}a(\tau)\Bigr]^2+\mathbb{P}(\tau)\Bigl[\frac{n\omega c^\theta\theta}{1-n}(c+\tau)^{-\theta-1}\Bigr]\\
=&\mathbb{P}(\tau)\Biggl[\Bigl(\omega+\frac{n\omega}{1-n}a(\tau)\Bigr)^2+\frac{n\omega}{1-n}\Phi(\tau)\Biggr]\\
=&\exp\Biggl\{-\omega\tau-\frac{n\omega}{1-n}A(\tau)\Biggr\}\Biggl[\Bigl(\omega+\frac{n\omega}{1-n}a(\tau)\Bigr)^2+\frac{n\omega}{1-n}\Phi(\tau)\Biggr].
\end{align*}
In conclusion,
\begin{align}
\label{eqn:sottocaso100}
F_{\tau}(\tau)=&\frac{1}{\lambda}\exp\Biggl\{-\omega\tau-\frac{n\omega}{1-n}A(\tau)\Biggr\}\Biggl[\Bigl(\omega+\frac{n\omega}{1-n}a(\tau)\Bigr)^2+\frac{n\omega}{1-n}\Phi(\tau)\Biggr].
\end{align}
We can see that this density doesn't depend on the assumed model for the conditional distribution of triggered events' magnitude. Therefore, in the case $m=m_0$, in which all the events are observable, our  hypothesis doesn't play any role.

\section{\textbf{The general case $m\ge m_0$}}
In order to find the expression of the inter-event time density for the general case, let's start again with formula~\eqref{eqn:b18}. Let's define the function
\begin{equation}
\label{eqn:psi}
\Psi(h(\cdot),\widetilde{m})=\int_{\widetilde{m}}^\infty p(m')e^{-\rho(m')h(m')}dm',
\end{equation}
where $\widetilde{m}\ge0$ and $h(\cdot)$ is a continuous function dependent on some variables; it is here considered as a function of only one of them, corresponding to the triggering event's magnitude $m'$. Recalling expression~\eqref{eqn:b13} of the probability generating function $G(z;t,\tau,m')$, we can rewrite the integral $\int_{m_0}^\infty G(z;t,\tau,m')p(m')dm'$ of the last member of~\eqref{eqn:b18} as
\begin{equation}
\label{eqn:uno}
\int_{m_0}^\infty G(z;t,\tau,m')p(m')dm'=\Psi[y(z;t,\tau,m_0,m'),m_0],
\end{equation}
where
\begin{align}
\label{eqn:y1}
y(z;t,\tau,m_0,m')&=b(t+\tau)-\int_0^tdx\,\Phi(x)\int_{m_0}^\infty G(z;t-x,\tau,m'')p(m''|m')dm''\notag\\
&-\int_0^{\tau}dx\,\Phi(t+\tau-x)\int_{m_0}^\infty G(z;0,x,m'')p(m''|m')dm''\notag\\
&+(1-z)\int_0^{\tau}dx\,\Phi(t+\tau-x)\int_{m}^\infty G(z;0,x,m'')p(m''|m')dm''\notag\\
=& b(t+\tau)-(\Phi(\cdot)\otimes D(z;\cdot,\tau,m_0,m'))(t)\notag\\
&-(\Phi(t+\cdot)\otimes D_+(z;\cdot,m_0,m'))(\tau)\notag\\
&+(1-z)(\Phi(t+\cdot)\otimes D_+(z;\cdot,m,m'))(\tau).
\end{align}
In the latter equality we have set
\begin{equation}
\label{eqn:b19a}
D(z;t,\tau,m_0,\overline{m})=\int_{m_0}^\infty G(z;t,\tau,m')p(m'|\overline{m})dm',
\end{equation}
\begin{equation}
\label{eqn:b19c}
D_+(z;t,\widetilde{m},\overline{m})=\int_{\widetilde{m}}^\infty G(z;0,t,m')p(m'|\overline{m})dm'.
\end{equation}
Let's notice that in our case we have $\widetilde{m}=m_0,m$.

Similarly, the other two integrals with respect to magnitude
\[
\int_{m_0}^\infty G(z;0,t,m')p(m')dm'\quad\text{and}\quad\int_m^\infty G(z;0,t,m')p(m')dm',
\]
in the last member of~\eqref{eqn:b18}, can be expressed as
\begin{equation}
\begin{split}
\label{eqn:due}
\int_{m_0}^\infty G(z;0,t,m')p(m')dm'&=\Psi[s(z;t,m_0,m'),m_0],\\
\int_{m}^\infty G(z;0,t,m')p(m')dm'&=\Psi[s(z;t,m_0,m'),m],
\end{split}
\end{equation}
where
\begin{align}
\label{eqn:s1}
s(z;t,m_0,m')=&b(t)-\int_0^tdx\,\Phi(t-x)\int_{m_0}^\infty G(z;0,x,m'')p(m''|m')dm''\notag\\
&+(1-z)\int_0^tdx\,\Phi(t-x)\int_{m}^\infty G(z;0,x,m'')p(m''|m')dm''\notag\\
=& b(t)-(\Phi(\cdot)\otimes D_+(z;\cdot,m_0,m'))(t)\notag\\
&+(1-z)(\Phi(\cdot)\otimes D_+(z;\cdot,m,m'))(t).
\end{align}
Concluding, equation~\eqref{eqn:b18} becomes
\begin{align}
\label{eqn:lll}
L(z;\tau)=&\int_0^\infty\Bigl\{1-\Psi[y(z;t,\tau,m_0,m'),m_0]\Bigr\}dt\notag\\
&+\int_0^\tau\Bigl\{1-\Psi[s(z;t,m_0,m'),m_0]+(1-z)\Psi[s(z;t,m_0,m'),m]\Bigr\}dt.
\end{align}

As explained before, in order to find the probability $\mathbb{P}(\tau)$ of having zero events in a temporal interval of length $\tau$, we have to compute the probability generating function in zero. We then have
\begin{equation}
\label{eqn:c1}
\mathbb{P}(\tau)=e^{-\omega L(0;\tau)}.
\end{equation}
Obviously, the expression of the function $L(0;\tau)$ can be obtained by setting $z=0$ in~\eqref{eqn:lll}. In particular we get
\begin{align}
\label{eqn:c2}
L(0;\tau)=&\int_0^\infty\Bigl\{1-\Psi[y(0;t,\tau,m_0,m'),m_0]\Bigr\}dt\notag\\
&+\int_0^\tau\Bigl\{1-\Psi[s(0;t,m_0,m'),m_0]+\Psi[s(0;t,m_0,m'),m]\Bigr\}dt\notag\\
=&\int_0^\infty\overline{N}(t,\tau)dt+\int_0^\tau\overline{N}_-(t)dt,
\end{align}
where
\begin{align}
\label{eqn:c4}
\overline{N}(t,\tau)=&1-\Psi[y(0;t,\tau,m_0,m'),m_0]
\end{align}
and
\begin{align}
\label{eqn:c55}
\overline{N}_-(t)=&1-\Psi[s(0;t,m_0,m'),m_0]+\Psi[s(0;t,m_0,m'),m].
\end{align}
\begin{observation}
Recalling equations~\eqref{eqn:uno} and~\eqref{eqn:due}, it is clear that we can rewrite in the following way the last two functions just defined:
\begin{align}
\label{eqn:altrenbarra1}
\overline{N}(t,\tau)&=1-\int_{m_0}^\infty G(0;t,\tau,m')p(m')dm'
\end{align}
and
\begin{align}
\label{eqn:altrenbarra2}
\overline{N}_-(t)&=1-\int_{m_0}^\infty G(0;0,t,m')p(m')dm'+\int_{m}^\infty G(0;0,t,m')p(m')dm'.
\end{align}
Obviously, in the case $m=m_0$ we have that $\overline{N}(t,\tau)=\overline{N}_{m=m_0}(t,\tau)$ and $\overline{N}_-(t)=1$.
\end{observation}
At this point, in order to obtain the probability $\mathbb{P}(\tau)$ of having zero events in $[0,\tau]$, we have to find an expression for the functions $\overline{N}(t,\tau)$ and $\overline{N}_-(t)$ and then substitute the obtained results in equation~\eqref{eqn:c2}. Actually, since they are integrals rather difficult to compute, we can approximate them by means of the first order Taylor power series expansion of the function $\Psi(h(\cdot),\widetilde{m})$, with respect to the first argument, around the point zero. In our case it holds $\widetilde{m}=m_0,m$ and, recalling equations~\eqref{eqn:y1} and~\eqref{eqn:s1}, the function $h(\cdot)$ in $\overline{N}(t,\tau)$ and $\overline{N}_-(t)$ is respectively equal to $y(0;t,\tau,m_0,\cdot)$ and $s(0;t,m_0,\cdot)$, where
\begin{align}
\label{eqn:y}
y(0;t,\tau,m_0,m')=&b(t+\tau)-(\Phi(\cdot)\otimes D(0;\cdot,\tau,m_0,m'))(t)\notag\\
&-(\Phi(t+\cdot)\otimes D_+(0;\cdot,m_0,m'))(\tau)\notag\\
&+(\Phi(t+\cdot)\otimes D_+(0;\cdot,m,m'))(\tau)
\end{align}
and
\begin{align}
\label{eqn:s}
s(0;t,m_0,m')=&b(t)-(\Phi(\cdot)\otimes D_+(0;\cdot,m_0,m'))(t)\notag\\
&+(\Phi(\cdot)\otimes D_+(0;\cdot,m,m'))(t).
\end{align}
We notice that the previous approximation makes sense since, in both cases, the first argument of $\Psi$ is small when $\tau\approx0$. In fact, we obtain that
\begin{itemize}
\item for the argument $y(0;t,\tau,m_0,m')$ holds
      \begin{align*}
      y(0;t,\tau=0,m_0,m')=&b(t)-b(t)=0.
      \end{align*}
      This follows because, recalling that $G(z;t,\tau,m'')$ is the probability generating function of the number of events triggered in $[0,\tau]$ by a shock with magnitude $m'$, occurred $t$ seconds before time 0, we have
      \begin{align*}
      D(0;t,\tau=0,m_0,m')=&\int_{m_0}^\infty G(0;t,\tau=0,m'')p(m''|m')dm''\\
      =&\int_{m_0}^\infty\bigl[0^0\mathbb{P}(N(0)=0)\\
      &+0^1\mathbb{P}(N(0)=1)\bigr]p(m''|m')dm''\\
      =&1.
      \end{align*}
      The previous calculus is valid because of the regularity of the process and the fact that the probability $\mathbb{P}(N(0)=0)$ of having zero events in an interval of length zero is equal to one.
\item on the other hand, for $s(0;t,m_0,m')$ holds
      \begin{align*}
      s(0;t=0,m_0,m')=&0,
      \end{align*}
      because $t$ varies between 0 and $\tau=0$.
\end{itemize}
Let's compute then the above-cited expansion and the relative approximations of $\overline{N}(t,\tau)$ and $\overline{N}_-(t)$. Let's start with the first function:
\begin{align}
\label{eqn:aap}
\overline{N}(t,\tau)=&1-\Psi[y(0;t,\tau,m_0,m'),m_0]\notag\\
=&1-\int_{m_0}^\infty p(m')e^{-\rho(m')y(0;t,\tau,m_0,m')}dm'\notag\\
\approx&1-\Biggl[1-\int_{m_0}^\infty p(m')\rho(m')y(0;t,\tau,m_0,m')dm'\Biggr]\notag\\
=&\int_{m_0}^\infty p(m')\rho(m')y(0;t,\tau,m_0,m')dm'\notag\\
=&\int_{m_0}^\infty p(m')\rho(m')dm'b(t+\tau)\notag\\
&-\int_{m_0}^\infty p(m')\rho(m')(\Phi(\cdot)\otimes D(0;\cdot,\tau,m_0,m'))(t)dm'\notag\\
&-\int_{m_0}^\infty p(m')\rho(m')(\Phi(t+\cdot)\otimes D_+(0;\cdot,m_0,m'))(\tau)dm'\notag\\
&+\int_{m_0}^\infty p(m')\rho(m')(\Phi(t+\cdot)\otimes D_+(0;\cdot,m,m'))(\tau)dm'\notag,\\
\end{align}
where we have used equations~\eqref{eqn:psi},~\eqref{eqn:c4} and~\eqref{eqn:y}.
We are going to study separately the integrals in the last equality of~\eqref{eqn:aap}. In the matter of the first, recalling the first equality in~\eqref{eqn:enne}, by definition we have
\[
\int_{m_0}^\infty p(m')\rho(m')dm'b(t+\tau)=nb(t+\tau).
\]
For the second one, recalling~\eqref{eqn:Agrande} and~\eqref{eqn:b19a}, we have
\begin{align*}
&\int_{m_0}^\infty p(m')\rho(m')(\Phi(\cdot)\otimes D(0;\cdot,\tau,m_0,m'))(t)dm'\\
&=\int_{m_0}^\infty dm'p(m')\rho(m')\int_0^tdx\,\Phi(t-x)\int_{m_0}^\infty G(0;x,\tau,m'')p(m''|m')dm''\\
&=\int_0^tdx\,\Phi(t-x)\int_{m_0}^\infty\int_{m_0}^\infty p(m')\rho(m')G(0;x,\tau,m'')p(m''|m')dm''dm'\\
&=\int_0^tdx\,\Phi(t-x)\int_{m_0}^\infty dm''G(0;x,\tau,m'')\int_{m_0}^\infty p(m')\rho(m')p(m''|m')dm'\\
&=\int_0^tdx\,\Phi(t-x)\int_{m_0}^\infty G(0;x,\tau,m'')I_A(m'')dm''.
\end{align*}
If we furthermore use equation~\eqref{eqn:b19c}, we can compute the third integral of the last expression in~\eqref{eqn:aap} as follows:
\begin{align*}
&\int_{m_0}^\infty p(m')\rho(m')(\Phi(t+\cdot)\otimes D_+(0;\cdot,m_0,m'))(\tau)dm'\\
&=\int_{m_0}^\infty dm'p(m')\rho(m')\int_0^{\tau}dx\,\Phi(t+\tau-x)\int_{m_0}^\infty G(0;0,x,m'')p(m''|m')dm''\\
&=\int_0^{\tau}dx\,\Phi(t+\tau-x)\int_{m_0}^\infty\int_{m_0}^\infty  p(m')\rho(m')G(0;0,x,m'')p(m''|m')dm''dm'\\
&=\int_0^{\tau}dx\,\Phi(t+\tau-x)\int_{m_0}^\infty dm''G(0;0,x,m'')\int_{m_0}^\infty p(m')\rho(m')p(m''|m')dm'\\
&=\int_0^{\tau}dx\,\Phi(t+\tau-x)\int_{m_0}^\infty G(0;0,x,m'')I_A(m'')dm''.
\end{align*}
Similarly, the fourth integral becomes
\begin{align*}
&\int_{m_0}^\infty p(m')\rho(m')(\Phi(t+\cdot)\otimes D_+(0;\cdot,m,m'))(\tau)dm'\\
&=\int_{m_0}^\infty dm'p(m')\rho(m')\int_0^{\tau}dx\,\Phi(t+\tau-x)\int_{m}^\infty G(0;0,x,m'')p(m''|m')dm''\\
&=\int_0^{\tau}dx\,\Phi(t+\tau-x)\int_{m_0}^\infty\int_{m}^\infty  p(m')\rho(m')G(0;0,x,m'')p(m''|m')dm''dm'\\
&=\int_0^{\tau}dx\,\Phi(t+\tau-x)\int_{m}^\infty dm''G(0;0,x,m'')\int_{m_0}^\infty p(m')\rho(m')p(m''|m')dm'\\
&=\int_0^{\tau}dx\,\Phi(t+\tau-x)\int_{m}^\infty G(0;0,x,m'')I_A(m'')dm''.
\end{align*}
We can then find the approximated expression of the function $\overline{N}(t,\tau)$ by substituting in equation~\eqref{eqn:aap} the obtained results. We get
\begin{align}
\label{eqn:approx1}
\overline{N}(t,\tau)\approx&nb(t+\tau)-\int_0^{t}dx\,\Phi(t-x)\int_{m_0}^\infty G(0;x,\tau,m'')I_A(m'')dm''\notag\\
&-\int_0^{\tau}dx\,\Phi(t+\tau-x)\int_{m_0}^\infty G(0;0,x,m'')I_A(m'')dm''\notag\\
&+\int_0^{\tau}dx\,\Phi(t+\tau-x)\int_{m}^\infty G(0;0,x,m'')I_A(m'')dm''.
\end{align}
For simplicity of notation, from now on we will substitute approximation symbol with equal sign.

Before proceeding with the study of the approximation of the second function $\overline{N}_-(t)$ in~\eqref{eqn:c2}, we focus for a moment on the approximation obtained for $\overline{N}(t,\tau)$. As we're going to see, it can be demonstrated that
\begin{align}
\label{eqn:mm02}
\int_0^\infty\overline{N}(t,\tau)dt&=\frac{n}{1-n}(a(\cdot)\otimes\overline{N}_-(\cdot))(\tau).
\end{align}
Consequently, the function $L(0;\cdot)$ in~\eqref{eqn:c2} can be written in terms of the function $\overline{N}_-(\cdot)$ only. We have indeed
\begin{align}
\label{eqn:mm03}
L(0;\tau)&=\int_0^\infty\overline{N}(t,\tau)dt+\int_0^\tau\overline{N}_-(t)dt\notag\\
&=\frac{n}{1-n}(a(\cdot)\otimes\overline{N}_-(\cdot))(\tau)+\int_0^\tau\overline{N}_-(t)dt.
\end{align}
In order to demonstrate equation~\eqref{eqn:mm02}, we impose in~\eqref{eqn:approx1} that is $I_A(m'')=np(m'')$, as already assumed in our model. This condition, as explained in Observation~\ref{oss:invperpcond}, guarantees the validity of the Gutenberg-Richter law when averaging on all triggering events' magnitudes. We then get
\begin{align}
\overline{N}(t,\tau)=&nb(t+\tau)-n\int_0^t dx\,\Phi(x)\int_{m_0}^\infty G(0;t-x,\tau,m'')p(m'')dm''\notag\\
&-n\int_0^\tau dx\,\Phi(t+\tau-x)\int_{m_0}^\infty G(0;0,x,m'')p(m'')dm''\notag\\
&+n\int_0^\tau dx\,\Phi(t+\tau-x)\int_{m}^\infty G(0;0,x,m'')p(m'')dm''\notag\\
=&n[b(t+\tau)-b(t)]+n\biggl\{\int_0^tdx\,\Phi(x)\biggl[1-\int_{m_0}^\infty G(0;t-x,\tau,m'')p(m'')dm''\biggr]\biggr\}\notag\\
&+n\biggl\{\int_0^\tau dx\,\Phi(t+\tau-x)\biggl[-\int_{m_0}^\infty G(0;0,x,m'')p(m'')dm''\notag\\
&+\int_{m}^\infty G(0;0,x,m'')p(m'')dm''\biggr]\biggr\}\notag\\
=&n\int_{t}^{t+\tau}\Phi(y)dy+n\int_0^t\Phi(x)\overline{N}(t-x,\tau)dx\notag\\
&+n\biggl\{\int_{t}^{t+\tau}dy\,\Phi(y)\biggl[-\int_{m_0}^\infty G(0;0,t+\tau-y,m'')p(m'')dm''\notag\\
&+\int_{m}^\infty G(0;0,t+\tau-y,m'')p(m'')dm''\biggr]\biggr\}\notag\\
=&n(\Phi(\cdot)\otimes\overline{N}(\cdot,\tau))(t)+n\int_{t}^{t+\tau}dy\,\Phi(y)\bigg[1\notag\\
&-\int_{m_0}^\infty G(0;0,t+\tau-y,m'')p(m'')dm''\notag\\
&+\int_{m}^\infty G(0;0,t+\tau-y,m'')p(m'')dm''\biggr]\notag\\
=&n(\Phi(\cdot)\otimes\overline{N}(\cdot,\tau))(t)+n\int_{t}^{t+\tau}\Phi(y)\overline{N}_-(t+\tau-y)dy\notag\\
=&n(\Phi(\cdot)\otimes\overline{N}(\cdot,\tau))(t)+n\int_{0}^{\tau}\Phi(t+\tau-x)\overline{N}_-(x)dx\notag,
\end{align}
where we have used equations~\eqref{eqn:altrenbarra1} and~\eqref{eqn:altrenbarra2} and in addition, in third and last equality, we have made the change of variable $y=t+\tau-x$. It then follows that
\begin{align}
\label{eqn:mm01}
\overline{N}(t,\tau)=&n(\Phi(\cdot)\otimes\overline{N}(\cdot,\tau))(t)+n(\Phi(t+\cdot)\otimes\overline{N}_-(\cdot))(\tau).
\end{align}
Integrating both the members with respect to times, between zero and infinity, we get
\begin{align}
\int_0^\infty \overline{N}(t,\tau)dt=&n\int_0^\infty\int_0^t\Phi(t-x)\overline{N}(x,\tau)dxdt\notag\\
&+n\int_0^\infty\int_0^{\tau}\Phi(t+\tau-x)\overline{N}_-(x)dxdt\notag\\
=&n\int_0^\infty\int_x^{\infty}\Phi(t-x)\overline{N}(x,\tau)dtdx\notag\\
&+n\int_0^{\tau}\int_0^\infty\Phi(t+\tau-x)\overline{N}_-(x)dtdx\notag\\
=&n\int_0^\infty dx\overline{N}(x,\tau)\int_x^{\infty}\Phi(t-x)dt\notag\\
&+n\int_0^{\tau}dx\overline{N}_-(x)\int_0^\infty\Phi(t+\tau-x)dt\notag\\
=&n\int_0^\infty\overline{N}(x,\tau)dx+n\int_0^{\tau}dx\overline{N}_-(x)\int_{\tau-x}^\infty\Phi(y)dy\notag\\
=&n\int_0^\infty\overline{N}(x,\tau)dx+n\int_0^{\tau}\overline{N}_-(x)a(\tau-x)dx\notag\\
=&n\int_0^\infty\overline{N}(x,\tau)dx+n(a(\cdot)\otimes\overline{N}_-(\cdot))(\tau),\notag
\end{align}
where we have used equations~\eqref{eqn:omdens} and~\eqref{eqn:funa}. We have then obtained equation~\eqref{eqn:mm02}.

At this point, we are ready to proceed with the study of the approximation of $\overline{N}_-(t)$. Using equations~\eqref{eqn:psi},~\eqref{eqn:c55} and~\eqref{eqn:s} we get
\begin{align}
\label{eqn:app}
\overline{N}_-(t)=&1-\Psi[s(0;t,m_0,m'),m_0]+\Psi[s(0;t,m_0,m'),m]\notag\\
=&1-\int_{m_0}^\infty p(m')e^{-\rho(m')s(0;t,m_0,m')}dm'+\int_{m}^\infty p(m')e^{-\rho(m')s(0;t,m_0,m')}dm'\notag\\
\approx&1-\Biggl[1-\int_{m_0}^\infty p(m')\rho(m')s(0;t,m_0,m')dm'\Biggr]\notag\\
&+\Biggl[\int_{m}^\infty p(m')dm'-\int_{m}^\infty p(m')\rho(m')s(0;t,m_0,m')dm'\Biggr]\notag\\
=&\int_{m_0}^\infty p(m')\rho(m')s(0;t,m_0,m')dm'+Q\notag\\
&-\int_{m}^\infty p(m')\rho(m')s(0;t,m_0,m')dm'\notag\\
=&Q+\int_{m_0}^\infty p(m')\rho(m')dm'b(t)\notag\\
&-\int_{m_0}^\infty p(m')\rho(m')(\Phi(\cdot)\otimes D_+(0;\cdot,m_0,m'))(t)dm'\notag\\
&+\int_{m_0}^\infty p(m')\rho(m')(\Phi(\cdot)\otimes D_+(0;\cdot,m,m'))(t)dm'\notag\\
&-\int_{m}^\infty p(m')\rho(m')dm'b(t)\notag\\
&+\int_{m}^\infty p(m')\rho(m')(\Phi(\cdot)\otimes D_+(0;\cdot,m_0,m'))(t)dm'\notag\\
&-\int_{m}^\infty p(m')\rho(m')(\Phi(\cdot)\otimes D_+(0;\cdot,m,m'))(t)dm',
\end{align}
where
\begin{equation}
\label{eqn:Qgrande}
Q=\int_{m}^\infty p(m')dm'=e^{-\beta(m-m_0)}.
\end{equation}
As previously done, we study separately all the integrals in the final expression of ~\eqref{eqn:app}. Regarding the first three integrals, recalling again equations~\eqref{eqn:Agrande},~\eqref{eqn:enne} and~\eqref{eqn:b19c}, we obtain the following results:
\[
\int_{m_0}^\infty p(m')\rho(m')dm'b(t)=nb(t);
\]
\begin{align*}
&\int_{m_0}^\infty p(m')\rho(m')(\Phi(\cdot)\otimes D_+(0;\cdot,m_0,m'))(t)dm'\\
&=\int_{m_0}^\infty dm'p(m')\rho(m')\int_{0}^{t}dx\,\Phi(t-x)\int_{m_0}^\infty G(0;0,x,m'')p(m''|m')dm''\\
&=\int_{0}^{t}dx\,\Phi(t-x)\int_{m_0}^\infty dm''G(0;0,x,m'')\int_{m_0}^\infty p(m')\rho(m')p(m''|m')dm'\\
&=\int_{0}^{t}dx\,\Phi(t-x)\int_{m_0}^\infty G(0;0,x,m'')I_A(m'')dm'';
\end{align*}
and in the end, similarly,
\begin{align*}
\int_{m_0}^\infty &p(m')\rho(m')(\Phi(\cdot)\otimes D_+(0;\cdot,m,m'))(t)dm'=\\
&=\int_{0}^{t}dx\,\Phi(t-x)\int_{m}^\infty G(0;0,x,m'')I_A(m'')dm''.
\end{align*}
In the same way, for the last three integrals in~\eqref{eqn:app} we get the following expressions. The first one becomes
\begin{align*}
\int_{m}^\infty p(m')\rho(m')dm'b(t)&=b(t)\int_{m}^\infty\beta e^{-\beta(m'-m_0)}\kappa e^{a(m'-m_0)}dm'\notag\\
&=b(t)\beta\kappa\int_{m}^\infty e^{-(\beta-a)(m'-m_0)}dm'\notag\\
&=b(t)\frac{\beta\kappa}{-(\beta-a)} e^{-(\beta-a)(m'-m_0)}\bigg|^\infty_{m}\notag\\
&=\frac{\beta\kappa e^{-(\beta-a)(m-m_0)}}{\beta-a}b(t)\notag\\
&=ne^{-(\beta-a)(m-m_0)}b(t),
\end{align*}
where we have used again equation~\eqref{eqn:enne} and the ones~\eqref{eqn:GR} and~\eqref{eqn:produttivita} in Appendix~\ref{sec:appc}. The second integral becomes
\begin{align*}
&\int_{m}^\infty p(m')\rho(m')(\Phi(\cdot)\otimes D_+(0;\cdot,m_0,m'))(t)dm'\\
&=\int_{m}^\infty dm'p(m')\rho(m')\int_{0}^{t}dx\,\Phi(t-x)\int_{m_0}^\infty G(0;0,x,m'')p(m''|m')dm''\\
&=\int_{0}^{t}dx\,\Phi(t-x)\int_{m_0}^\infty dm''G(0;0,x,m'')\int_{m}^\infty p(m')\rho(m')p(m''|m')dm'\\
&=\int_{0}^{t}dx\,\Phi(t-x)\int_{m_0}^\infty G(0;0,x,m'')I_B(m'')dm'',
\end{align*}
where we have set
\begin{equation}
\label{eqn:Bgrande}
I_B(m'')=\int_{m}^\infty p(m')\rho(m')p(m''|m')dm'.
\end{equation}
In the end, the third one is similarly
\begin{align*}
\int_{m}^\infty &p(m')\rho(m')(\Phi(\cdot)\otimes D_+(0;\cdot,m,m'))(t)dm'=\\
&=\int_{0}^{t}dx\,\Phi(t-x)\int_{m}^\infty G(0;0,x,m'')I_B(m'')dm''.
\end{align*}
Substituting the obtained results in~\eqref{eqn:app}, we are then able to find the approximated expression of the function $\overline{N}_-(t)$:
\begin{align}
\label{eqn:approx2}
\overline{N}_-(t)\approx&Q+nb(t)\notag\\
&-\int_{0}^{t}dx\,\Phi(t-x)\int_{m_0}^\infty G(0;0,x,m'')I_A(m'')dm''\notag\\
&+\int_{0}^{t}dx\,\Phi(t-x)\int_{m}^\infty G(0;0,x,m'')I_A(m'')dm''\notag\\
&-ne^{-(\beta-a)(m-m_0)}b(t)\notag\\
&+\int_{0}^{t}dx\,\Phi(t-x)\int_{m_0}^\infty G(0;0,x,m'')I_B(m'')dm''\notag\\
&-\int_{0}^{t}dx\,\Phi(t-x)\int_{m}^\infty G(0;0,x,m'')I_B(m'')dm''\notag\\
=&Q+nb(t)\biggl[1-e^{-(\beta-a)(m-m_0)}\biggr]\notag\\
&-\int_{0}^{t}dx\,\Phi(t-x)\int_{m_0}^\infty G(0;0,x,m'')I_A(m'')dm''\notag\\
&+\int_{0}^{t}dx\,\Phi(t-x)\int_{m}^\infty G(0;0,x,m'')I_A(m'')dm''\notag\\
&+\int_{0}^{t}dx\,\Phi(t-x)\int_{m_0}^\infty G(0;0,x,m'')I_B(m'')dm''\notag\\
&-\int_{0}^{t}dx\,\Phi(t-x)\int_{m}^\infty G(0;0,x,m'')I_B(m'')dm''.
\end{align}

\begin{figure}[t]
\centering
\includegraphics[width=20pc]{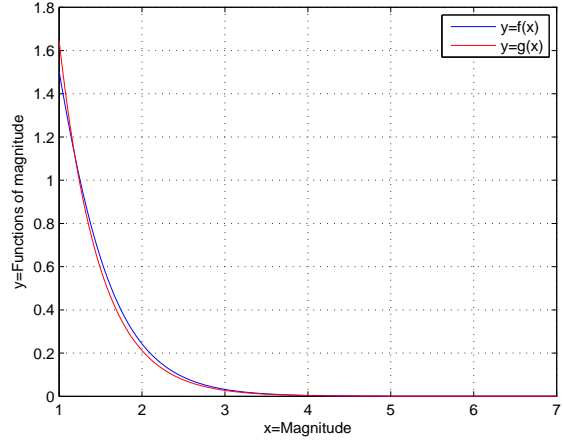}
\caption{Plot of the functions $f(x)=Hp(x)+C_1p(x)(1-2\frac{p(x)}{\beta})H(1-H)$ and $g(x)=Lp(x)$, where $p(x)=\beta\exp\{-\beta(x-m_0)\}$.}
\label{fig:Bint}
\end{figure}

At this point we can work on the approximation of the function $\overline{N}_-(t)$ as we made for $\overline{N}(t,\tau)$. In particular, the first condition to impose is $I_A(m'')=np(m'')$, justified in the Observation~\ref{oss:invperpcond}. This is guaranteed by our proposal of triggered events' magnitude conditional density function:
\begin{align*}
p(m''|m')&=p(m'')\bigl[1+f(m',m'')\bigr]\notag\\
&=\beta e^{-\beta(m''-m_0)}\biggl[1+C_1\bigl(1-2e^{-(\beta-a)(m'-m_0)}\bigr)\bigl(1-2e^{-\beta(m''-m_0)}\bigr)\biggr];
\end{align*}
for further details about this conditional density law see Appendix~\ref{sec:appa}. On the other hand, in order to find the condition that the integral $I_B(m'')$ must satisfy, we substitute in its expression, found in~\eqref{eqn:Bgrande}, the above-mentioned proposal. It holds
\begin{align}
\label{eqn:intb}
I_B(m'')=&\int_m^\infty p(m')\rho(m')p(m''|m')dm'\notag\\
=&p(m'')\int_m^\infty p(m')\rho(m')dm'+p(m'')\int_m^\infty p(m')\rho(m')f(m',m'')dm'\notag\\
=&p(m'')ne^{-(\beta-a)(m-m_0)}\notag\\
&+p(m'')\int_m^\infty p(m')\rho(m')C_1\bigl(1-2e^{-(\beta-a)(m'-m_0)}\bigr)\bigl(1-2e^{-\beta(m''-m_0)}\bigr)dm'\notag\\
=&p(m'')ne^{-(\beta-a)(m-m_0)}\notag\\
&+p(m'')C_1\bigl(1-2e^{-\beta(m''-m_0)}\bigr)\int_m^\infty p(m')\rho(m')\bigl(1-2e^{-(\beta-a)(m'-m_0)}\bigr)dm'\notag\\
=&p(m'')ne^{-(\beta-a)(m-m_0)}+p(m'')C_1\bigl(1-2e^{-\beta(m''-m_0)}\bigr)ne^{-(\beta-a)(m-m_0)}\notag\\
&-p(m'')C_1\bigl(1-2e^{-\beta(m''-m_0)}\bigr)2\int_m^\infty\beta\kappa e^{-2(\beta-a)(m'-m_0)}dm'\notag\\
=&p(m'')ne^{-(\beta-a)(m-m_0)}+p(m'')C_1\bigl(1-2e^{-\beta(m''-m_0)}\bigr)ne^{-(\beta-a)(m-m_0)}\notag\\
&+p(m'')C_1\bigl(1-2e^{-\beta(m''-m_0)}\bigr)\frac{\beta\kappa}{\beta-a}\int_m^\infty -2(\beta-a)e^{-2(\beta-a)(m'-m_0)}dm'\notag\\
=&p(m'')ne^{-(\beta-a)(m-m_0)}+p(m'')C_1\bigl(1-2e^{-\beta(m''-m_0)}\bigr)ne^{-(\beta-a)(m-m_0)}\notag\\
&+p(m'')C_1\bigl(1-2e^{-\beta(m''-m_0)}\bigr)n e^{-2(\beta-a)(m'-m_0)}\big|_m^\infty\notag\\
=&p(m'')ne^{-(\beta-a)(m-m_0)}+p(m'')C_1\bigl(1-2e^{-\beta(m''-m_0)}\bigr)ne^{-(\beta-a)(m-m_0)}\notag\\
&-p(m'')C_1\bigl(1-2e^{-\beta(m''-m_0)}\bigr)ne^{-2(\beta-a)(m-m_0)},\notag\\
=&p(m'')nH+p(m'')C_1\biggl(1-2\frac{p(m'')}{\beta}\biggr)nH-p(m'')C_1\biggl(1-2\frac{p(m'')}{\beta}\biggr)nH^2,
\end{align}
where the last equality follows observing that $e^{-\beta(m''-m_0)}=\frac{p(m'')}{\beta}$ and setting
\begin{align}
\label{eqn:mm04}
H=e^{-(\beta-a)(m-m_0)}.
\end{align}
The expression obtained for integral $I_B(m'')$ depends then both on the Gutenberg-Richter function $p(m'')$ and on its square. This complicates a lot our study. However, observing that $|1-2e^{-\beta(m''-m_0)}|<1$, $C_1<1$ and that the difference $H-H^2$ is generally small for the typical values of parameters, we can approximate $I_B(m'')$ with a function of kind $p(m'')nL$, where the constant $L$ is very close to $H$. As an example, we consider the following parameters data-set:
\begin{align*}
\beta&=2.0493,\\
a&=1.832,\\
C_1&=0.8,\\
m_0&=1,\\
m&=1.8.
\end{align*}
By means of a numerical algorithm minimizing the difference between $I_B(m'')$ and $p(m'')nL$, we get the constant $L=0.8032$ that is, as we expected, close to the constant $H=0.8404$. A graphical check is given in Fig.~\ref{fig:Bint}.

We can then proceed with the substitution in equation~\eqref{eqn:approx2} of the two conditions
\begin{align}
I_A(m'')&=np(m'')\notag\\
\label{eqn:mm07}
I_B(m'')&\approx nLp(m''),
\end{align}
for a certain constant $L$ to be found numerically.

Substituting again the symbol of approximation with equal sign for simplicity of notation, we have
\begin{align}
\overline{N}_-(t)=&Q+nb(t)[1-H]\notag\\
&-n\int_{0}^{t}dx\,\Phi(t-x)\int_{m_0}^\infty G(0;0,x,m'')p(m'')dm''\notag\\
&+n\int_{0}^{t}dx\,\Phi(t-x)\int_{m}^\infty G(0;0,x,m'')p(m'')dm''\notag\\
&+nL\int_{0}^{t}dx\,\Phi(t-x)\int_{m_0}^\infty G(0;0,x,m'')p(m'')dm''\notag\\
&-nL\int_{0}^{t}dx\,\Phi(t-x)\int_{m}^\infty G(0;0,x,m'')p(m'')dm''\notag\\
=&Q+nb(t)-nb(t)H+nLb(t)-nLb(t)\notag\\
&-n\int_{0}^{t}dx\,\Phi(x)\int_{m_0}^\infty G(0;0,t-x,m'')p(m'')dm''\notag\\
&+n\int_{0}^{t}dx\,\Phi(x)\int_{m}^\infty G(0;0,t-x,m'')p(m'')dm''\notag\\
&+nL\int_{0}^{t}dx\,\Phi(x)\int_{m_0}^\infty G(0;0,t-x,m'')p(m'')dm''\notag\\
&-nL\int_{0}^{t}dx\,\Phi(x)\int_{m}^\infty G(0;0,t-x,m'')p(m'')dm''\notag\\
=&Q-nb(t)H+nLb(t)\notag\\
&+n[1-L]\int_0^tdx\,\Phi(x)\biggl[1-\int_{m_0}^\infty G(0;0,t-x,m'')p(m'')dm''\notag\\
&+\int_{m}^\infty G(0;0,t-x,m'')p(m'')dm''\biggr]\notag\\
=&Q-nb(t)[H-L]+n\bigl[1-L\bigr]\int_0^t\Phi(x)\overline{N}_-(t-x)dx\notag,
\end{align}
where we have used equation~\eqref{eqn:altrenbarra2}. We then obtain
\begin{align}
\label{eqn:mm08}
\overline{N}_-(t)=&Q-nb(t)[H-L]+\delta(\Phi(\cdot)\otimes\overline{N}_-(\cdot))(t),
\end{align}
where
\begin{align}
\label{eqn:mm09}
\delta=n\bigl[1-L\bigr].
\end{align}

Before proceeding with the research of the function $\overline{N}_-(t)$, let's do the following observation.
\begin{observation}
If we integrate both the members of~\eqref{eqn:mm08} with respect to time between zero and infinity, recalling equations~\eqref{eqn:funa} and~\eqref{eqn:funagrande} we get
\begin{align*}
\int_0^{\tau}\overline{N}_-(t)dt=&Q\tau-n[H-L]\int_0^{\tau}b(t)dt+\delta\int_0^{\tau}\int_0^t\Phi(t-x)\overline{N}_-(x)dxdt\\
=&Q\tau-n[H-L]\int_0^{\tau}\Bigl[1-a(t)\Bigr]dt+\delta\int_0^{\tau}\int_x^{\tau}\Phi(t-x)\overline{N}_-(x)dtdx\\
=&Q\tau-n[H-L]\biggl[\tau-\int_0^{\tau}a(t)dt\biggr]+\delta\int_0^{\tau}dx\overline{N}_-(x)\int_0^{\tau-x}\Phi(y)dy\\
=&Q\tau-n[H-L]\biggl[\tau-A(\tau)\biggr]+\delta\int_0^{\tau}\overline{N}_-(x)[1-a(\tau-x)]dx\\
=&Q\tau-n[H-L]\biggl[\tau-A(\tau)\biggr]+\delta\int_0^{\tau}\overline{N}_-(x)dx-\delta(a(\cdot)\otimes\overline{N}_-(\cdot))(\tau).
\end{align*}
Hence, it is possible to find the convolution $(a(\cdot)\otimes\overline{N}_-(\cdot))(\tau)$ as a function of the integral $\int_0^{\tau}\overline{N}_-(t)dt$ and vice versa. In particular, we find
\begin{align}
\label{eqn:mm014}
\int_0^{\tau}\overline{N}_-(t)dt=&\frac{1}{1-\delta}\biggl\{\tau\Bigl[Q-n[H-L]\Bigr]+n[H-L]A(\tau)-\delta(a(\cdot)\otimes\overline{N}_-(\cdot))(\tau)\biggr\}
\end{align}
and
\begin{align}
\label{eqn:mm015}
(a(\cdot)\otimes\overline{N}_-(\cdot))(\tau)=&\frac{\tau}{\delta}\Bigl[Q-n[H-L]\Bigr]+\frac{n[H-L]}{\delta}A(\tau)+\Bigl[1-\frac{1}{\delta}\Bigr]\int_0^{\tau}\overline{N}_-(t)dt.
\end{align}
Using equation~\eqref{eqn:mm015}, we can rewrite the function $L(0;\tau)$ (equation~\eqref{eqn:mm03}) only in terms of $\int_0^{\tau}\overline{N}_-(t)dt$:
\begin{align}
\label{eqn:mm0101}
L(0;\tau)=&\frac{n}{1-n}\Biggl[\tau\frac{Q-n[H-L]}{\delta}+\frac{n[H-L]}{\delta}A(\tau)+\frac{\delta-1}{\delta}\int_0^{\tau}\overline{N}_-(t)dt\Biggl]+\int_0^\tau\overline{N}_-(t)dt\notag\\
=&\int_0^{\tau}\overline{N}_-(t)dt\Bigl[\frac{n}{1-n}\frac{\delta-1}{\delta}+1\Bigl]+\tau n\frac{Q-n[H-L]}{\delta(1-n)}+\frac{n^2[H-L]}{\delta(1-n)}A(\tau).
\end{align}
Recalling~\eqref{eqn:c1}, let's multiply equation~\eqref{eqn:mm0101} by the rate $\omega$ relative to background events. We obtain
\begin{align}
\label{eqn:mm0100}
\omega L(0;\tau)=&\frac{\delta-1}{\delta}\omega\Delta\int_0^{\tau}\overline{N}_-(t)dt+\tau\frac{n}{1-n}\frac{\widetilde{\omega}}{\delta}+\omega\frac{n^2[H-L]}{\delta(1-n)}A(\tau),
\end{align}
where we have set
\begin{align}
\label{eqn:deltagrande}
\Delta&=\frac{n}{1-n}-\frac{\delta}{1-\delta}\\
&\qquad\quad\text{and}\notag\\
\label{eqn:omegam}
\widetilde{\omega}&=\omega\bigl\{Q-n[H-L]\bigr\}.
\end{align}
\end{observation}

If we could have an expression for $\int_0^{\tau}\overline{N}_-(t)dt$, using equation~\eqref{eqn:mm0100} in~\eqref{eqn:c1} it would be possible to find the probability of having zero events in $[0,\tau]$.

In order to obtain an explicit expression of the function $\overline{N}_-(t)$, we use the standard technique of Laplace transform. More precisely, we integrate between zero and infinity both the members of equation~\eqref{eqn:mm08} multiplied by the suitable exponential function. We get
\begin{align*}
\int_0^\infty e^{-st}\overline{N}_-(t)dt=&\int_0^\infty e^{-st}Qdt-n[H-L]\int_0^\infty e^{-st}b(t)dt\\
&+\delta\int_0^\infty e^{-st}(\Phi(\cdot)\otimes\overline{N}_-(\cdot))(t)dt\\
=&\frac{Q}{s}-n[H-L]\int_0^\infty e^{-st}\Bigl[1-a(t)\Bigr]dt\\
&+\delta\int_0^\infty e^{-st}\Phi(t)dt\int_0^\infty e^{-st}\overline{N}_-(t)dt,\quad s\in\mathbb{C}.
\end{align*}
It then follows that
\begin{align}
\int_0^\infty e^{-st}\overline{N}_-(t)dt=&\frac{1}{1-\delta\int_0^\infty e^{-st}\Phi(t)dt}\biggl[\frac{Q-n[H-L]}{s}+n[H-L]\int_0^\infty e^{-st}a(t)dt\biggr].
\end{align}
Since
\begin{align*}
\int_0^\infty e^{-st}a(t)dt&=\int_0^\infty e^{-st}\int_t^\infty\Phi(x)dxdt\\
&=\int_0^\infty e^{-st}\biggl[1-\int_0^t\Phi(x)dx\biggr]dt\\
&=\int_0^\infty e^{-st}dt-\int_0^\infty e^{-st}\int_0^t\Phi(x)dxdt\\
&=-\frac{1}{s}e^{-st}\mid^\infty_0-\int_0^\infty\Phi(x)\int_x^{\infty}e^{-st}dtdx\\
&=\frac{1}{s}-\frac{1}{s}\int_0^\infty\Phi(x)e^{-sx}dx\\
&=\frac{1}{s}\biggl[1-\int_0^\infty e^{-sx}\Phi(x)dx\biggr],
\end{align*}
%
%
we obtain
\begin{align}
\label{eqn:mm010}
\int_0^\infty e^{-st}\overline{N}_-(t)dt=&\frac{1}{s[1-\delta\int_0^\infty e^{-st}\Phi(t)dt]}\biggl[Q-n[H-L]\int_0^\infty e^{-st}\Phi(t)dt\biggr].
\end{align}
We notice that the above function of the complex variable $s$ has one pole in $s=0$. Furthermore, there are no poles with real part larger than zero. In fact, the other factor in the denominator is never zero for $\Re(s)>0$, since $\delta<1$ and
\begin{align*}
\int_0^\infty e^{-st}\Phi(t)dt&<\bigg|\int_0^\infty e^{-st}\Phi(t)dt\bigg|\\
&<\int_0^\infty |e^{-st}\Phi(t)|dt\\
&<\int_0^\infty |\Phi(t)|dt\\
&=1.
\end{align*}

Let's then compute the Laplace transform of the function $\Phi(\cdot)$ appearing in~\eqref{eqn:mm010}. It is
\begin{align}
\label{eqn:mm012}
\int_0^\infty e^{-st}\Phi(t)dt&=\theta c^{\theta}\int_0^\infty e^{-st}(t+c)^{-1-\theta}dt\notag\\
&=\theta c^{\theta}\int_0^\infty e^{-x}\biggl(\frac{x}{s}+c\biggr)^{-1-\theta}\frac{1}{s}dx\notag\\
&=\theta (sc)^{\theta}e^{sc}\int_0^\infty e^{-(x+sc)}(x+sc)^{-\theta-1}dx\notag\\
&=\theta (sc)^{\theta}e^{sc}\int_{sc}^\infty e^{-y}y^{-\theta-1}dy\notag\\
&=\theta (sc)^{\theta}e^{sc}\Gamma(-\theta,sc),
\end{align}
where
\[
\Gamma(l,x)=\int_x^\infty e^{-t}t^{l-1}dt,\quad|\arg x|<\pi
\]
is the incomplete gamma function (see~\cite{temme:primo}). We notice that the integral in the third line of~\eqref{eqn:mm012} is the explicit form of the complex integral along the horizontal half-line starting from the complex point $sc$, written in the consecutive line.

Substituting in equation~\eqref{eqn:mm010} the solution~\eqref{eqn:mm012}, we obtain
\begin{align}
\label{eqn:fin}
\mathfrak{L}[\overline{N}_-(t)](s)&=\int_0^\infty e^{-st}\overline{N}_-(t)dt\notag\\
&=\frac{1}{s[1-\delta\theta (sc)^{\theta}e^{sc}\Gamma(-\theta,sc)]}\biggl[Q-n[H-L]\theta (sc)^{\theta}e^{sc}\Gamma(-\theta,sc)\biggr].
\end{align}
One possible way to find the function $\overline{N}_-(t)$ consists of computing the inverse Laplace transform through the Bromwich integral:
\begin{equation*}
\overline{N}_-(t)=\frac{1}{2\pi\imath}\int_{\sigma-\imath\infty}^{\sigma+\imath\infty}e^{st}\mathfrak{L}[\overline{N}_-(t)](s)ds,
\end{equation*}
where $\sigma$ is any positive number. We then get
\begin{align}
\label{eqn:br2}
\overline{N}_-(t)=&\frac{1}{2\pi}\int_{-\infty}^{+\infty}e^{(\sigma+\imath\xi)t}\mathfrak{L}[\overline{N}_-(t)](\sigma+\imath\xi)d\xi\notag\\
=&\frac{e^{\sigma t}}{2\pi}\int_{-\infty}^{+\infty}e^{\imath\xi t}\mathfrak{L}[\overline{N}_-(t)](\sigma+\imath\xi)d\xi\notag\\
=&\frac{e^{\sigma t}}{2\pi}\int_{-\infty}^{+\infty}\bigl[\mathfrak{L}[\overline{N}_-(t)](\sigma+\imath\xi)\cos\xi t+\imath\mathfrak{L}[\overline{N}_-(t)](\sigma+\imath\xi)\sin\xi t\bigr]d\xi\notag\\
=&\frac{e^{\sigma t}}{2\pi}\biggl\{\int_{-\infty}^{+\infty}\Bigl[\Re\bigl(\mathfrak{L}[\overline{N}_-(t)](\sigma+\imath\xi)\bigr)+\imath\Im\bigl(\mathfrak{L}[\overline{N}_-(t)](\sigma+\imath\xi)\bigr)\Bigr]\cos\xi td\xi\notag\\
&+\imath\int_{-\infty}^{+\infty}\Bigl[\Re\bigl(\mathfrak{L}[\overline{N}_-(t)](\sigma+\imath\xi)\bigr)+\imath\Im\bigl(\mathfrak{L}[\overline{N}_-(t)](\sigma+\imath\xi)\bigr)\Bigr]\sin\xi td\xi\Biggr\}\notag\\
=&\frac{e^{\sigma t}}{2\pi}\biggl\{\int_{-\infty}^{+\infty}\Bigl[\Re\bigl(\mathfrak{L}[\overline{N}_-(t)](\sigma+\imath\xi)\bigr)\cos\xi t-\Im\bigl(\mathfrak{L}[\overline{N}_-(t)](\sigma+\imath\xi)\bigr)\sin\xi t\Bigr]d\xi\notag\\
&+\imath\int_{-\infty}^{+\infty}\Bigl[\Im\bigl(\mathfrak{L}[\overline{N}_-(t)](\sigma+\imath\xi)\bigr)\cos\xi t+\Re\bigl(\mathfrak{L}[\overline{N}_-(t)](\sigma+\imath\xi)\bigr)\sin\xi t\Bigr]d\xi\biggr\},
\end{align}
where
\begin{align*}
\Re\bigl(\mathfrak{L}[\overline{N}_-(t)](\sigma+\imath\xi)\bigr)&=\Re\Biggl(\int_0^\infty\overline{N}_-(t)e^{-(\sigma+\imath\xi)t}dt\Biggr)=\int_0^\infty e^{-\sigma t}\overline{N}_-(t)\cos\xi tdt\\
&\text{and}\\
\Im\bigl(\mathfrak{L}[\overline{N}_-(t)](\sigma+\imath\xi)\bigr)&=\Im\Biggl(\int_0^\infty\overline{N}_-(t)e^{-(\sigma+\imath\xi)t}dt\Biggr)=-\int_0^\infty e^{-\sigma t}\overline{N}_-(t)\sin\xi tdt.
\end{align*}
Since $\overline{N}_-(t)$ is a real function, we have
\begin{equation*}
\int_{-\infty}^{+\infty}\Bigl[\Im\bigl(\mathfrak{L}[\overline{N}_-(t)](\sigma+\imath\xi)\bigr)\cos\xi t+\Re\bigl(\mathfrak{L}[\overline{N}_-(t)](\sigma+\imath\xi)\bigr)\sin\xi t\Bigr]d\xi=0;
\end{equation*}
actually, this also follows easily from the odd-parity of the integrand. Moreover, observing that the function
\[
W(\xi)=\Re\bigl(\mathfrak{L}[\overline{N}_-(t)](\sigma+\imath\xi)\bigr)\cos\xi t-\Im\bigl(\mathfrak{L}[\overline{N}_-(t)](\sigma+\imath\xi)\bigr)\sin\xi t
\]
is of even-parity, the last expression in~\eqref{eqn:br2} becomes:
\[
\frac{e^{\sigma t}}{\pi}\int_0^\infty\Bigl[\Re\bigl(\mathfrak{L}[\overline{N}_-(t)](\sigma+\imath\xi)\bigr)\cos\xi t-\Im\bigl(\mathfrak{L}[\overline{N}_-(t)](\sigma+\imath\xi)\bigr)\sin\xi t\Bigr]d\xi.
\]
Finally, recalling that $\overline{N}_-(t)=0$ for $t<0$, one can verify that
\begin{equation*}
\int_{-\infty}^{+\infty}\Bigl[\Re\bigl(\mathfrak{L}[\overline{N}_-(t)](\sigma+\imath\xi)\bigr)\cos\xi t+\Im\bigl(\mathfrak{L}[\overline{N}_-(t)](\sigma+\imath\xi)\bigr)\sin\xi t\Bigr]d\xi=0.
\end{equation*}
Then,
\begin{align}
\label{eqn:br3}
\overline{N}_-(t)&=\frac{2e^{\sigma t}}{\pi}\int_0^\infty\Re\bigl(\mathfrak{L}[\overline{N}_-(t)](\sigma+\imath\xi)\bigr)\cos\xi td\xi\notag\\
&=\frac{2e^{\sigma t}}{\pi}\int_0^\infty\Re\Biggl(\frac{Q-n[H-L]\theta [(\sigma+\imath\xi)c]^{\theta}e^{(\sigma+\imath\xi)c}\Gamma(-\theta,(\sigma+\imath\xi)c)}{(\sigma+\imath\xi)[1-\delta\theta [(\sigma+\imath\xi)c]^{\theta}e^{(\sigma+\imath\xi)c}\Gamma(-\theta,(\sigma+\imath\xi)c)]}\Biggr)\cos\xi td\xi.
\end{align}
One may use a numerical approach to compute the last integral and then to find the required function $\overline{N}_-(t)$. Once done this, we must integrat with respect to time. Then, the obtained result must be substituted in equation~\eqref{eqn:mm0100}. In conclusion, inter-event time density is obtained recalling that
\begin{equation}
\label{eqn:h}
F_{\tau}(\tau)=\frac{1}{\lambda}\frac{d^2}{d\tau^2}e^{-\omega L(0;\tau)}.
\end{equation}

Because of the difficulty of both the computation of $\overline{N}_-(t)$ and the inter-event time density, we simplify the case we are analyzing by a further approximation.

If we consider again $\tau$ sufficiently small, we can set $a(\tau-t)\approx a(\tau)$ in~\eqref{eqn:mm014}. Substituting once again the symbol of approximation with the one of equality for simplicity of notation, we get
\begin{align*}
(1-\delta)\int_0^{\tau}\overline{N}_-(t)dt=&\tau\Bigl[Q-n[H-L]\Bigr]+n[H-L]A(\tau)-\delta\int_0^{\tau}a(\tau)\overline{N}_-(t)dt,
\end{align*}
from which we have
\begin{align}
\label{eqn:mm016}
\int_0^{\tau}\overline{N}_-(t)dt=&\frac{1}{1-\delta+\delta a(\tau)}\biggl\{\tau\Bigl[Q-n[H-L]\Bigr]+n[H-L]A(\tau)\biggr\}.
\end{align}
The above-mentioned approximation makes sense since, if $t\in[0,\tau]$ and $\tau$ in small, so is $t$.

We can then conveniently substitute the expression~\eqref{eqn:mm016} of $\int_0^{\tau}\overline{N}_-(t)dt$ in~\eqref{eqn:mm0100}:
\begin{align}
\label{eqn:mm019}
\omega L(0;\tau)=&\frac{\delta-1}{\delta}\omega\Delta\frac{1}{1-\delta+\delta a(\tau)}\biggl\{\tau\Bigl[Q-n[H-L]\Bigr]+n[H-L]A(\tau)\biggr\}\notag\\
&+\tau\widetilde{\omega}\frac{1}{\delta}\frac{n}{1-n}+\omega\frac{n}{\delta(1-n)} n[H-L]A(\tau)\notag\\
=&\frac{\delta-1}{\delta}\widetilde{\omega}\Delta\tau\frac{1}{1-\delta+\delta a(\tau)}\notag\\
&+\frac{\delta-1}{\delta}\omega\Delta\frac{1}{1-\delta+\delta a(\tau)}n[H-L]A(\tau)\notag\\
&+\tau\widetilde{\omega}\frac{1}{\delta}\frac{n}{1-n}+\omega\frac{n}{\delta(1-n)}n[H-L]A(\tau)\notag\\
=&\frac{\tau\widetilde{\omega}}{\delta}\biggl[\frac{(\delta-1)\Delta}{1-\delta+\delta a(\tau)}+\frac{n}{1-n}\biggr]\notag\\
&+\frac{\omega n[H-L]A(\tau)}{\delta}\biggl[\frac{(\delta-1)\Delta}{1-\delta+\delta a(\tau)}+\frac{n}{1-n}\biggr]\notag\\
=&\biggl[\tau\widetilde{\omega}+\omega n[H-L]A(\tau)\biggr]\biggl[\frac{\frac{\delta-1}{\delta}\Delta}{1-\delta+\delta a(\tau)}+\frac{n}{(1-n)\delta}\biggr]\notag\\
=&\biggl[\tau\widetilde{\omega}+\omega n[H-L]A(\tau)\biggr]\biggl[\frac{1-n+na(\tau)}{(1-n)[1-\delta+\delta a(\tau)]}\biggr],
\end{align}
where we have used equation~\eqref{eqn:omegam} and the last equality follows because
\begin{align*}
\frac{\frac{\delta-1}{\delta}\Delta}{1-\delta+\delta a(\tau)}+\frac{n}{(1-n)\delta}=&\frac{\delta-1}{\delta}\biggl[\frac{n}{1-n}+\frac{\delta}{\delta-1}\biggr]\frac{1}{1-\delta+\delta a(\tau)}+\frac{n}{(1-n)\delta}\\
=&\biggl[\frac{n(\delta-1)}{(1-n)\delta}+1\biggr]\frac{1}{1-\delta+\delta a(\tau)}+\frac{n}{(1-n)\delta}\\
=&\frac{n(\delta-1)+(1-n)\delta}{(1-n)\delta[1-\delta+\delta a(\tau)]}+\frac{n}{(1-n)\delta}\\
=&\frac{n(\delta-1)+(1-n)\delta+n(1-\delta)+n\delta a(\tau)}{(1-n)\delta[1-\delta+\delta a(\tau)]}\\
=&\frac{1-n+na(\tau)}{(1-n)[1-\delta+\delta a(\tau)]}.
\end{align*}
In this last calculus we have used equation~\eqref{eqn:deltagrande}.

Now, using~\eqref{eqn:h}, we are able to find the probability density function $F_{\tau}(\tau)$ relative to inter-event time. From the calculations to obtain the second derivative of $\mathbb{P}(\tau)=e^{-\omega L(0;\tau)}$ (equation~\eqref{eqn:derivsec}), shown in Appendix~\ref{sec:apph}, we get
\begin{align}
\label{eqn:mm020}
F_{\tau}(\tau)=&\frac{1}{\lambda}e^{-\omega L(0;\tau)}\Biggl\{\biggl[\biggl(\widetilde{\omega}+\omega n[H-L]a(\tau)\biggr)\biggl(\frac{1-n+na(\tau)}{(1-n)[1-\delta+\delta a(\tau)]}\biggr)\notag\\
&+\biggl(\tau\widetilde{\omega}+\omega n[H-L]A(\tau)\biggr)\biggl(\frac{\Phi(\tau)(\delta-n)}{(1-n)[1-\delta+\delta a(\tau)]^2}\biggr)\biggr]^2\notag\\
&-\biggl[-\omega n[H-L]\Phi(\tau)\biggl(\frac{1-n+na(\tau)}{(1-n)[1-\delta+\delta a(\tau)]}\biggr)\notag\\
&+2\biggl(\widetilde{\omega}+\omega n[H-L]a(\tau)\biggr)\biggl(\frac{\Phi(\tau)(\delta-n)}{(1-n)[1-\delta+\delta a(\tau)]^2}\biggr)\notag\\
&+\biggl(\tau\widetilde{\omega}+\omega n[H-L]A(\tau)\biggr)\frac{(\delta-n)\Phi(\tau)[(\theta+1)(\delta-1)+\delta a(\tau)(\theta-1)]}{(1-n)(\tau+c)[1-\delta+\delta a(\tau)]^3}\biggr]\Biggr\}\notag\\
=&\frac{1}{\lambda}\exp\Biggl\{-\biggl[\tau\widetilde{\omega}+\omega n[H-L]A(\tau)\biggr]\biggl[\frac{1-n+na(\tau)}{(1-n)[1-\delta+\delta a(\tau)]}\biggr]\Biggr\}\notag\\
&\cdot\Biggl\{\biggl[\biggl(\widetilde{\omega}+\omega n[H-L]a(\tau)\biggr)\biggl(\frac{1-n+na(\tau)}{(1-n)[1-\delta+\delta a(\tau)]}\biggr)\notag\\
&+\biggl(\tau\widetilde{\omega}+\omega n[H-L]A(\tau)\biggr)\biggl(\frac{\Phi(\tau)(\delta-n)}{(1-n)[1-\delta+\delta a(\tau)]^2}\biggr)\biggr]^2\notag\\
&-\biggl[-\omega n[H-L]\Phi(\tau)\biggl(\frac{1-n+na(\tau)}{(1-n)[1-\delta+\delta a(\tau)]}\biggr)\notag\\
&+2\biggl(\widetilde{\omega}+\omega n[H-L]a(\tau)\biggr)\biggl(\frac{\Phi(\tau)(\delta-n)}{(1-n)[1-\delta+\delta a(\tau)]^2}\biggr)\notag\\
&+\biggl(\tau\widetilde{\omega}+\omega n[H-L]A(\tau)\biggr)\frac{(\delta-n)\Phi(\tau)[(\theta+1)(\delta-1)+\delta a(\tau)(\theta-1)]}{(1-n)(\tau+c)[1-\delta+\delta a(\tau)]^3}\biggr]\Biggr\}.
\end{align}
We notice that, differently from the case $m=m_0$, in this case $m>m_0$, in which not all the events are observable, our hypothesis of a conditional probability density function for the magnitudes of triggered events plays a role. In fact, this conditional distribution appears in expression~\eqref{eqn:mm020} through some of the constants, as for example $L$.

\section{\textbf{Conclusions and future work}}
We proposed and analyzed a new version of the ETAS model, where triggered event's magnitudes are probabilistically dependent on triggering event's ones. Thanks to the tool of the probability generating function and the Palm theory, we obtained a closed-form for the density of inter-event time for small values. More precisely, the closed-form is found by a linear approximation around zero of an exponential function. This is the only approximation needed for the case $m=m_0$. In the general case $m\ge m_0$, we have used an additional approximation to get the result in a closed-form. However, not using the second approximation, we were able to find the studied density function in terms of the inverse Laplace transform of a suitable function. The results obtained show that our hypothesis of conditioning plays a role when we do not observe all the triggering events, as happens in practice. In fact, if we consider the particular case $m=m_0$, that is if we observe all the events that are able to trigger, the density obtained for the inter-event doesn't depend on the new proposal for the distribution of triggered events' magnitude. It depends only on the Omori-Utsu law and on its temporal integral. On the other hand, in the case $m>m_0$, in which not all the triggering events are observable, we can see the influence of the above-mentioned proposal on some constants appearing in $F_{\tau}(\tau)$. The inter-event time density then depends on the Omori-Utsu law, on its temporal integrals and on the new probability density function for triggered events' magnitude. A possible problem to be analyzed for future work is the asymptotic study of the inter-event density for long times. Furthermore, it could be interesting to implement the new model and then to perform a simulation study. The statistical inference issue could be also of interest.

\section*{\textbf{Appendix}}

\appendix

\section{\textbf{Phenomenological laws}}
\label{sec:appc}
The \emph{Gutenberg-Richter law}
\begin{align}
\label{eqn:GR}
p(m')&=b\ln10\cdot10^{-b(m'-m_0)}\notag\\
&=b\ln10\cdot e^{-b\ln10(m'-m_0)}\notag\\
&=\beta e^{-\beta(m'-m_0)},
\end{align}
with $\beta=b\ln10$, describes the empirical distribution of magnitudes. The value $m_0$ is the minimum magnitude an event must have to be able to trigger its own progeny. This law is adopted for events' magnitude when we do not consider the respective characteristics of past seismicity.

The \emph{Omori-Utsu law}
\begin{equation}
\label{eqn:omori}
\Phi(t)=\frac{\theta c^{\theta}}{(c+t)^{1+\theta}}
\end{equation}
can be interpreted as the probability density function of random times at which first generation shocks independently occur, when the triggering event has occurred in $t=0$.

The \emph{productivity law}
\begin{align}
\label{eqn:produttivita}
\rho(m')&=\kappa10^{\alpha(m'-m_0)}\notag\\
&=\kappa e^{\alpha\ln10(m'-m_0)}\notag\\
&=\kappa e^{a(m'-m_0)},
\end{align}
where $a=\alpha\ln10$, represents the contribution of first generation shocks triggered by an event with magnitude equals to $m'$.

For the previous three laws, the parameters $(\beta,\theta,c,\kappa,a)$ are positive and typically estimated using maximum-likelihood algorithms.

\section{\textbf{Probability generating function and some of its fundamental properties}}
\label{sec:appb}
The probability generating function is a very useful tool to study random variables. It has several properties among which we can enunciate the following ones, which are very important for our analysis.
\begin{proposition}
Let $\{X_1,X_2,\dots\}$ be a sequence of independent and identically distributed random variables, with generating function $G_X$. Let also $N$ be a discrete positive random variable independent of $X_i\quad\forall i$, with generating function $G_N$. Then, given the random variable $S$ defined as
\[
S=\sum_{i=1}^N X_i,
\]
one has
\begin{equation}
\label{prop:1}
G_S(z)=G_N[G_X(z)].
\end{equation}
\end{proposition}
\begin{proposition}
Given two independent random variables $X$ and $Y$, it holds
\begin{equation}
\label{prop:2}
G_{X+Y}(z)=G_X(z)G_Y(z).
\end{equation}
\end{proposition}

\section{\textbf{Conditional probability density function for triggered events' magnitude}}
\label{sec:appa}
The probability density function we propose for triggered events' magnitude is not the Gutenberg-Richter law. Indeed, it is more realistic to suppose that the law of a triggered event's magnitude depends on the one of the corresponding triggering event. We then consider a conditional probability density function with respect to triggering event's magnitude. In what follows the latter magnitude will be indicated with $m'$. We search for a function $p(\cdot|m')$ such that:
\begin{align}
\label{eqn:siadensita}
\int_{m_0}^\infty &p(m''|m')dm''=1;\\
\notag\\
\label{eqn:condizioneperro}
&p(m''|tr)=p(m''),
\end{align}
where
\begin{equation}
\label{eqn:pperfigli}
p(m''|tr):=\int_{m_0}^\infty \frac{p(m')\rho(m')p(m''|m')}{\int_{m_0}^\infty p(m')\rho(m')dm'}dm'
\end{equation}
indicates the probability for an event to have magnitude $m''$, given the fact that it is not spontaneous. Recall that $p(m')$ is the Gutenberg-Richter law (equation~\eqref{eqn:GR} in Appendix~\ref{sec:appc}), $\rho(m')$ is the productivity law (equation~\eqref{eqn:produttivita} in Appendix~\ref{sec:appc}) and $m_0$ is the reference magnitude, that is the minimum value for an event to trigger other shocks. The above-cited condition~\eqref{eqn:condizioneperro} corresponds to the need of obtaining the Gutenberg-Richter law when we average over all the triggering event's magnitudes, as previously explained in Observation~\ref{oss:invperpcond}.
\begin{observation}
\label{oss:3}
The right side of definition~\eqref{eqn:pperfigli} can be derived as below. Let's consider the three following random variables:
\begin{itemize}
\item $M'\sim$ Gutenberg-Richter law is the triggering events' magnitude;
\item $N$ counts the number of shocks that a generic event triggers; it is such that $p(N=k)=\int_{m_0}^\infty p(N=k|m')p(m')dm'$;
\item $N_{m_{\delta}''}$ counts how many triggered shocks among the previous $N$ have magnitude $m_{\delta}''\in[m'',m''+\delta)$.
\end{itemize}
Since
\[
N_{m_{\delta}''}|(m',N=k)\sim Bin_{N_{m_{\delta}''}}\biggl(k,\int_{m''}^{m''+\delta}p(x|m')dx\biggr),
\]
where $Bin_{X}(n,p)$ indicates the law of a binomial random variable $X$ with parameters $n$ and $p$, we have
\[
p(N_{m_{\delta}''}=q)=\int_{m_0}^\infty\sum_{k=0}^\infty Bin_{N_{m_{\delta}''}=q}\biggl(k,\int_{m''}^{m''+\delta}p(x|m')dx\biggr)p(N=k|m')p(m')dm'.
\]

Let's consider the probability for a triggered event to have magnitude in the interval $[m'',m''+\delta)$. In order to find an expression of the just-mentioned probability, we suppose to have $\widetilde{n}$ realizations of the three variables $(M',N,N_{m_{\delta}''})$.
We focus on the ratio between the sample means of the random variables $N^i_{m_{\delta}''}$ and $N^i$, for $i=1,\dots,\widetilde{n}$. On one side, this ratio converges by the law of large numbers to the above probability. On the other side, it also converges in probability to the ratio of the respective expected values, that are the quantities
\begin{align}
\label{eqn:perinv}
&\frac{\sum_{q=0}^\infty q\int_{m_0}^\infty\sum_{k=0}^\infty Bin_{N_{m_{\delta}''}=q}\biggl(k,\int_{m''}^{m''+\delta}p(x|m')dx\biggr)p(N=k|m')p(m')dm'}{\sum_{k=0}^\infty k\int_{m_0}^\infty p(N=k|m')p(m')dm'}\notag\\
&=\frac{\int_{m_0}^\infty\sum_{k=0}^\infty\sum_{q=0}^\infty q Bin_{N_{m_{\delta}''}=q}\biggl(k,\int_{m''}^{m''+\delta}p(x|m')dx\biggr)p(N=k|m')p(m')dm'}{\int_{m_0}^\infty\sum_{k=0}^\infty k p(N=k|m')p(m')dm'}\notag\\
&=\frac{\int_{m_0}^\infty\sum_{k=0}^\infty k\int_{m''}^{m''+\delta}p(x|m')dx\,p(N=k|m')p(m')dm'}{\int_{m_0}^\infty \rho(m')p(m')dm'}\notag\\
&=\frac{\int_{m_0}^\infty\int_{m''}^{m''+\delta}p(x|m')dx\sum_{k=0}^\infty kp(N=k|m')p(m')dm'}{\int_{m_0}^\infty \rho(m')p(m')dm'}\notag\\
&=\frac{\int_{m_0}^\infty\int_{m''}^{m''+\delta}p(x|m')dx\rho(m')p(m')dm'}{\int_{m_0}^\infty \rho(m')p(m')dm'}.
\end{align}
Dividing by $\delta$ and computing the limit for $\delta$ tending to zero, we get precisely the probability for a triggered event to have magnitude $m''$:
\[
p(m''|tr)=\frac{\int_{m_0}^\infty p(m''|m')\rho(m')p(m')dm'}{n},
\]
where we recall that
\[
n=\int_{m_0}^\infty p(m')\rho(m')dm'=\frac{\beta\kappa}{\beta-a}.
\]
In equation~\eqref{eqn:perinv} we have used the fact that the productivity law $\rho(m')$, defined in Appendix~\ref{sec:appc} (equation~\eqref{eqn:produttivita}), can be rewritten as
\[
\rho(m')=\sum_{k=0}^\infty k p(N=k|m').
\]
In fact, it represents the mean number of first generation shocks triggered by an event with magnitude $m'$.
\end{observation}
We want to focus now on condition~\eqref{eqn:condizioneperro}, which can be rewritten in the following way:
\begin{equation}
\label{eqn:star}
\int_{m_0}^\infty p(m')\rho(m')p(m''|m')dm'=np(m'').
\end{equation}
Given the expressions of the functions $p(m')$ and $\rho(m')$, let's find a suitable law $p(m''|m')$ such that~\eqref{eqn:star} holds. Some qualitative behavior of this function has been obtained from the results of a statistical analysis of some Italian catalogues. The function we are looking for should be such that, when $m'$ increases, the probability of having events with high magnitudes must increase and, at the same time, the one of events with lower magnitudes must decrease obviously. For further information, see~\cite{mio:primo}.

In view of these considerations, the choice we adopt for $p(m''|m')$ is
\begin{equation}
\label{eqn:probcond1}
p(m''|m')=p(m'')\bigl[1+f(m',m'')\bigr],
\end{equation}
with
\begin{equation}
\label{eqn:f1}
f(m',m'')=-q(m')+s(m')\Bigl(1-e^{-\beta(m''-m_0)}\Bigr).
\end{equation}
For now, the functions $q(m')$ and $s(m')$ are unknown. As we are going to see, they can be chosen in order to satisfy the following conditions:
\begin{itemize}
\item $\int_0^\infty p(m''|m')=1$,
\item $p(m''|m')\ge0$,
\item $p(m''|m')$ is such that condition~\eqref{eqn:star} is verified,
\item $p(m''|m')$ has the qualitative behavior previously introduced.
\end{itemize}
Let's proceed then with the search of the above-mentioned functions $q(m')$ and $s(m')$. First of all, the conditional function $p(m''|m')$ must be a density. We then have
\begin{equation}
\label{eqn:cond1}
\int_{m_0}^\infty p(m'')\bigl[1+f(m',m'')\bigr]dm''=1\quad\Leftrightarrow\quad\int_{m_0}^\infty p(m'')f(m',m'')dm''=0.
\end{equation}
Substituting the expression chosen for $f(m',m'')$ we obtain
\begin{align*}
0&=\int_{m_0}^\infty p(m'')\bigl[-q(m')+s(m')-s(m')e^{-\beta(m''-m_0)}\bigr]dm''\\
&=-q(m')+s(m')-s(m')\int_{m_0}^\infty\beta e^{-2\beta(m''-m_0)}dm''\\
&=-q(m')+s(m')-\frac{s(m')}{2}\\
&=-q(m')+\frac{s(m')}{2},
\end{align*}
from which
\[
s(m')=2q(m').
\]
Substituting the last result in equation~\eqref{eqn:f1} we get
\begin{align}
\label{eqn:f2}
f(m',m'')&=q(m')-2q(m')e^{-\beta(m''-m_0)}\notag\\
&=q(m')\bigl(1-2e^{-\beta(m''-m_0)}\bigr).
\end{align}
Recalling equation~\eqref{eqn:probcond1}, in order to have also $p(m''|m')\ge0$, it should hold $f(m',m'')\ge-1$, that is
\[
q(m')\bigl(2e^{-\beta(m''-m_0)}-1\bigr)\le1.
\]
Since $|2e^{-\beta(m''-m_0)}-1|\le1$, this is guaranteed if
\begin{equation}
\label{condperq}
|q(m')|\le1.
\end{equation}

Let's impose at this point condition~\eqref{eqn:star}:
\begin{align}
\label{eqn:cond2}
\int_{m_0}^\infty p(m')\rho(m')p(m'')&\bigl[1+f(m',m'')\bigr]dm'=np(m'')\notag\\
&\Updownarrow\notag\\
\int_{m_0}^\infty p(m')\rho(m')&f(m',m'')dm'=0.
\end{align}
For $f(m',m'')$ as in~\eqref{eqn:f2} we get
\begin{align*}
0&=\int_{m_0}^\infty p(m')\rho(m')q(m')\bigl(1-2e^{-\beta(m''-m_0)}\bigr)dm'\\
&=\bigl(1-2e^{-\beta(m''-m_0)}\bigr)\int_{m_0}^\infty p(m')\rho(m')q(m')dm'.
\end{align*}
A sufficient condition to solve the previous equation is obviously
\begin{equation}
\label{eqn:perq}
\int_{m_0}^\infty p(m')\rho(m')q(m')dm'=0.
\end{equation}
Before making a choice of the function $q(m')$ such that condition~\eqref{eqn:perq} holds, let's make some qualitative comments.

\begin{figure}[t]
\centering
\includegraphics[width=20pc]{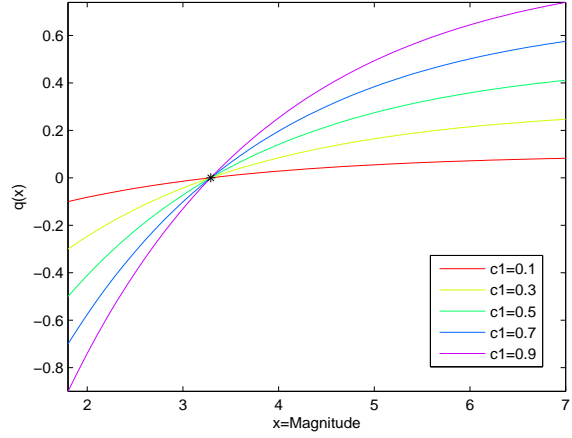}
\caption{Plot of the function $q(x)=C_1(1-2\exp\{-0.4648(x-1.8)\})$. The value at which the function becomes zero is $m^*=3.2913$.}
\label{fig:q}
\end{figure}

The function $q(m')$ should be such that the conditional law $p(m''|m')$ has the qualitative behavior before introduced. To this aim, we assume $q(m')$ to be continuous and increasing. Furthermore we impose that $q(m')$ is negative for $m'<\overline{m}$ at which it becomes zero and positive elsewhere (see Fig.~\ref{fig:q}).
Let's notice that when $m'=\overline{m}$ we have $f(m',m'')=0$ and the conditional probability density becomes the Gutenberg-Richter law.
From equation~\eqref{eqn:f2}, it can be observed that, since the function
\[
1-2e^{-\beta(m''-m_0)}
\]
is always increasing, the increasing or decreasing of $f(m',\cdot)$ will be determined by the sign of $q(m')$.
In particular, for $m'>\overline{m}$ where $q(m')>0$, the function $f(m',m'')$ increases in $m''$. For smaller (higher) values of
\[
m''=\frac{1}{\beta}\ln2+m_0,
\]
where $1-2e^{-\beta(m''-m_0)}=0$, the function $f(m',m'')$ will have negative (positive) sign. Hence, the conditional density $p(m''|m')=p(m'')[1+f(m',m'')]$ is below $p(m'')$ for magnitude values smaller than $m''=\frac{1}{\beta}\ln2+m_0$ and above it for higher magnitude values.
Vice versa, for $m'<\overline{m}$ where $q(m')<0$, the function $f(m',m'')$ is decreasing in $m''$. With analogous reasonings we come to the conclusion that $p(m''|m')$ is above the Gutenberg-Richter law for values smaller than $m''=\frac{1}{\beta}\ln2+m_0$, below for higher values.
\begin{observation}
\label{oss:4}
For triggered events' magnitude values such that $m''>\frac{1}{\beta}\ln2+m_0$, the factor $1-2e^{-\beta(m''-m_0)}$ is positive. Hence, fixed $m''$ larger than the just-mentioned value, since $q(m')$ is increasing, the function $f(m',m'')$ is increasing with respect to the first argument, too. This means that when triggering event's magnitude increases, so does the probability density function at a point $m''>\frac{1}{\beta}\ln2+m_0$. On the contrary, the value of the probability density decreases at a value $m''<\frac{1}{\beta}\ln2+m_0$. We are then in qualitative agreement with the graphical results obtained.
\end{observation}

To sum up, the function $q(m')$ must be:
\begin{itemize}
\item less than one in absolute value,
\item increasing,
\item negative for magnitudes less than a certain value $\overline{m}$, positive otherwise,
\item such that condition~\eqref{eqn:perq} is verified.
\end{itemize}
The choice we make for $q(m')$ is the following:
\begin{equation}
\label{eqn:q}
q(m')=-C_1+C_2\bigl(1-e^{-(\beta-a)(m'-m_0)}\bigr),
\end{equation}
where the constants $C_1$ and $C_2$ are obtained imposing the above-cited conditions for $q(\cdot)$. We specify that in practice $\beta>a$: it is the condition which ensures that the process doesn't explode.
Firstly, let's consider equation~\eqref{eqn:perq}. Recalling the expressions of $p(m')$ and $\rho(m')$ (respectively equations~\eqref{eqn:GR} and~\eqref{eqn:produttivita} in Appendix~\ref{sec:appc}), we get
\begin{align*}
0&=\beta\kappa\int_{m_0}^\infty e^{-(\beta-a)(m'-m_0)}\biggl[-C_1+C_2\bigl(1-e^{-(\beta-a)(m'-m_0)}\bigr)\biggr]dm'\\
&=\beta\kappa\biggl[(-C_1+C_2)\int_{m_0}^\infty e^{-(\beta-a)(m'-m_0)}dm'-C_2\int_{m_0}^\infty e^{-2(\beta-a)(m'-m_0)}dm'\biggr]\\
&=\beta\kappa\biggl[\frac{-C_1+C_2}{\beta-a}-\frac{C_2}{2(\beta-a)}\biggr]\\
&=\frac{\beta\kappa}{\beta-a}\biggl(-C_1+\frac{C_2}{2}\biggr),
\end{align*}
from which follows that
\[
C_2=2C_1.
\]
Hence
\begin{equation}
\label{eqn:q2}
q(m')=C_1-2C_1e^{-(\beta-a)(m'-m_0)}=C_1\bigl(1-2e^{-(\beta-a)(m'-m_0)}\bigr).
\end{equation}
In order to have also $|q(m')|<1$, since $|1-2e^{-(\beta-a)(m'-m_0)}|<1$, it is enough to choose $1>C_1\ge0$. One can easily verify that the function $q(\cdot)$ is increasing since $C_1\ge0$ and $\beta>a$ (Fig.~\ref{fig:q}). Furthermore, if $C_1\ne0$, the value for which $q(m')=0$ is
\[
m'=\overline{m}=\frac{1}{\beta-a}\ln2+m_0.
\]
Since $q(m')$ increases, it will be positive when $m'$ is larger than this value and is negative elsewhere.

At this point, substituting the expression~\eqref{eqn:q2} of $q(m')$ in~\eqref{eqn:f2}, we get
\begin{equation}
\label{eqn:f3}
f(m',m'')=C_1\bigl(1-2e^{-(\beta-a)(m'-m_0)}\bigr)\bigl(1-2e^{-\beta(m''-m_0)}\bigr).
\end{equation}
\begin{figure}[t]
\centering
\subfloat[][]
{\includegraphics[width=20pc]{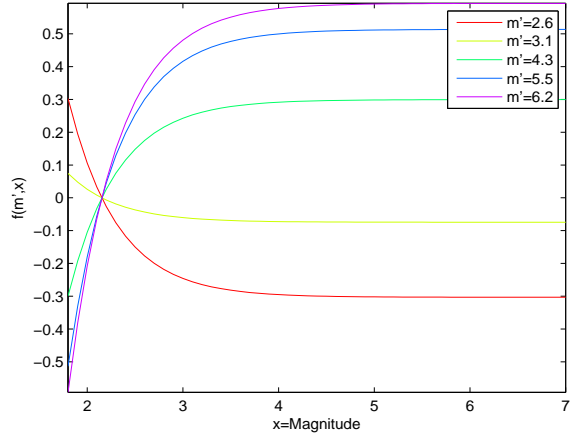}} \quad
\subfloat[][]
{\includegraphics[width=20pc]{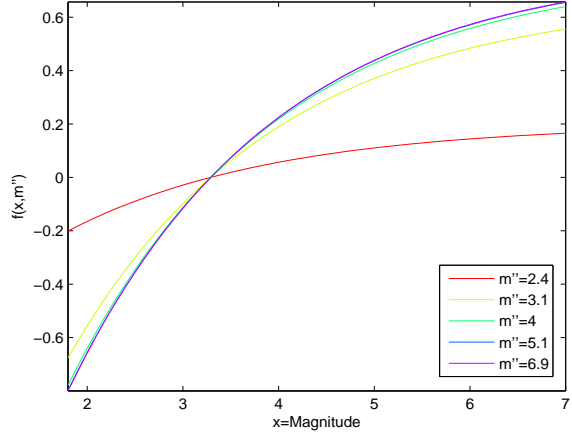}}
\caption{Plots of the functions $f(m',x)=0.8(1-2\exp\{-0.4648(m'-1.8)\})(1-2\exp\{-1.9648(x-1.8)\})$ and $f(x,m'')=0.8(1-2\exp\{-0.4648(x-1.8)\})(1-2\exp\{-1.9648(m''-1.8)\})$, respectively in (a) and (b). In the first subplot the value that varies is the triggering events' magnitude. On the other hand, in the second one the varying value is the triggered events' magnitude.}
\label{fig:funzionef}
\end{figure}
In Fig.~\ref{fig:funzionef}(a) we can see the plot of the function $f(m',\cdot)$ obtained with the values of parameters such that we get $\frac{\ln2}{\beta-a}+m_0=3.2913$. This is the value at which $q(m')$ becomes zero. The values of parameters correspond to a real situation because they were estimated based on a real data set. As expected, $f(m',m'')$ decreases (increases) with respect to $m''$ when $m'$ is smaller (higher) than the just-mentioned value. Furthermore, it holds $f(m',m'')=0$ for $m''=\frac{\ln2}{\beta}+m_0=2.1527$. In Fig.~\ref{fig:funzionef}(b) one can observe that, for any fixed value of $m''>2.1527$, the function $f(m',m'')$ is increasing with respect to the first variable.

In conclusion, substituting equations~\eqref{eqn:f3} and~\eqref{eqn:GR} (Appendix~\ref{sec:appc}) in~\eqref{eqn:probcond1}, we obtain the probability density function relative to triggered event's magnitude, conditioned on the magnitude $m'$ of its own triggering event:
\begin{align}
\label{eqn:pcond}
p(m''|m')&=\beta e^{-\beta(m''-m_0)}\biggl[1+C_1\bigl(1-2e^{-(\beta-a)(m'-m_0)}\bigr)\bigl(1-2e^{-\beta(m''-m_0)}\bigr)\biggr],
\end{align}
for all $1>C_1\ge0$.

\begin{figure}[h!]
\centering
\includegraphics[width=20pc]{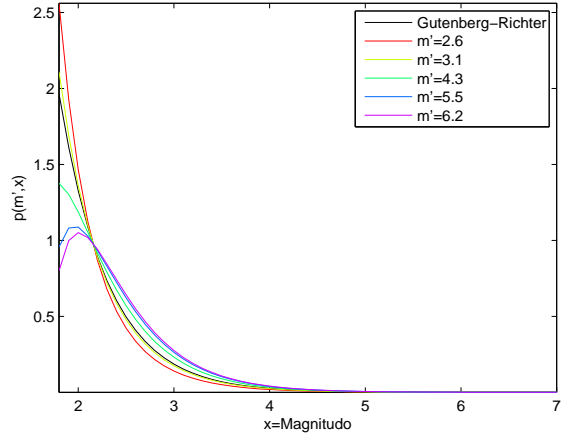}
\caption{Plot of the probability density function relative to triggered events' magnitude when triggering events' magnitude varies and $C_1=0.8$.}
\label{fig:pcond}
\end{figure}

We study now the behavior of $p(m''|m')$ by looking at its derivative. Recalling equations~\eqref{eqn:probcond1} and~\eqref{eqn:f2} we have
\begin{align*}
\frac{dp(m''|m')}{dm''}&=\frac{dp(m'')}{dm''}[1+f(m',m'')]+p(m'')\frac{df(m',m'')}{dm''}\\
&=-\beta p(m'')[1+f(m',m'')]+p(m'')\frac{df(m',m'')}{dm''}\\
&=\beta p(m'')\bigl[-1-f(m',m'')+2q(m')e^{-\beta(m''-m_0)}\bigr]\\
&=\beta p(m'')\bigl[-1-q(m')+2q(m')e^{-\beta(m''-m_0)}+2q(m')e^{-\beta(m''-m_0)}\bigr]\\
&=\beta p(m'')\bigl[-1-q(m')+4q(m')e^{-\beta(m''-m_0)}\bigr].
\end{align*}
It follows that if $q(m')\le0$, that is $m'\le\frac{1}{\beta-a}\ln2+m_0$, since $|q(m')|<1$, the density function $p(m''|m')$ is always decreasing in $m''$. If instead $q(m')>0$, that is $m'>\frac{1}{\beta-a}\ln2+m_0$, the above-cited density increases in $m''$ till a certain maximum reached in
\[
4q(m')e^{-\beta(m''-m_0)}=1+q(m')\quad\Leftrightarrow\quad m''=\frac{1}{\beta}\ln\frac{4q(m')}{1+q(m')}+m_0
\]
and after that it decreases. This kind of trend is shown in Fig.~\ref{fig:pcond}. This behavior agrees with the results of the analysis we made in~\cite{mio:primo}.
In Fig.~\ref{fig:pcond} we can also verify the behavior described in Observation~\ref{oss:4}.
\clearpage

\section{\textbf{Time derivatives of the probability $\mathbb{P}(\tau)$ of having zero events in $[0,\tau]$}}
\label{sec:apph}
In this appendix we include the calculations we made to obtain the second time derivative of the probability $\mathbb{P}(\tau)=e^{-\omega L(0;\tau)}$ in the case $m\ge m_0$. This derivative is indeed useful in order to get the inter-event time density $F_{\tau}(\tau)$ by means of equation~\eqref{eqn:intensintert}.
Using equations~\eqref{eqn:derivA} and~\eqref{eqn:mm019} we have
\begin{align*}
\frac{d}{d\tau}\omega L(0;\tau)=&\biggl[\widetilde{\omega}+\omega n[H-L]a(\tau)\biggr]\biggl[\frac{1-n+na(\tau)}{(1-n)[1-\delta+\delta a(\tau)]}\biggr]\\
&+\biggl[\tau\widetilde{\omega}+\omega n[H-L]A(\tau)\biggr]\frac{d}{d\tau}\biggl[\frac{1-n+na(\tau)}{(1-n)[1-\delta+\delta a(\tau)]}\biggr].
\end{align*}
It follows that
\begin{align*}
\frac{d\mathbb{P}(\tau)}{d\tau}=&-e^{-\omega L(0;\tau)}\Biggl\{\biggl[\widetilde{\omega}+\omega n[H-L]a(\tau)\biggr]\biggl[\frac{1-n+na(\tau)}{(1-n)[1-\delta+\delta a(\tau)]}\biggr]\\
&+\biggl[\tau\widetilde{\omega}+\omega n[H-L]A(\tau)\biggr]\frac{d}{d\tau}\biggl[\frac{1-n+na(\tau)}{(1-n)[1-\delta+\delta a(\tau)]}\biggr]\Biggr\}.
\end{align*}
The second derivative of $\mathbb{P}(\tau)$ with respect to $\tau$ is then
\begin{align}
\label{eqn:derivsec}
\frac{d^2\mathbb{P}(\tau)}{d\tau^2}=&e^{-\omega L(0;\tau)}\frac{d}{d\tau}\omega L(0;\tau)\Biggl\{\biggl[\widetilde{\omega}+\omega n[H-L]a(\tau)\biggr]\biggl[\frac{1-n+na(\tau)}{(1-n)[1-\delta+\delta a(\tau)]}\biggr]\notag\\
&+\biggl[\tau\widetilde{\omega}+\omega n[H-L]A(\tau)\biggr]\frac{d}{d\tau}\biggl[\frac{1-n+na(\tau)}{(1-n)[1-\delta+\delta a(\tau)]}\biggr]\Biggr\}\notag\\
&-e^{-\omega L(0;\tau)}\Biggl\{-\omega n[H-L]\Phi(\tau)\biggl[\frac{1-n+na(\tau)}{(1-n)[1-\delta+\delta a(\tau)]}\biggr]\notag\\
&+\biggl[\widetilde{\omega}+\omega n[H-L]a(\tau)\biggr]\frac{d}{d\tau}\biggl[\frac{1-n+na(\tau)}{(1-n)[1-\delta+\delta a(\tau)]}\biggr]\notag\\
&+\biggl[\widetilde{\omega}+\omega n[H-L]a(\tau)\biggr]\frac{d}{d\tau}\biggl[\frac{1-n+na(\tau)}{(1-n)[1-\delta+\delta a(\tau)]}\biggr]\notag\\
&+\biggl[\tau\widetilde{\omega}+\omega n[H-L]A(\tau)\biggr]\frac{d^2}{d\tau^2}\biggl[\frac{1-n+na(\tau)}{(1-n)[1-\delta+\delta a(\tau)]}\biggr]\Biggr\}\notag\\
=&e^{-\omega L(0;\tau)}\Biggl\{\biggl[\biggl(\widetilde{\omega}+\omega n[H-L]a(\tau)\biggr)\biggl(\frac{1-n+na(\tau)}{(1-n)[1-\delta+\delta a(\tau)]}\biggr)\notag\\
&+\biggl(\tau\widetilde{\omega}+\omega n[H-L]A(\tau)\biggr)\biggl(\frac{\Phi(\tau)(\delta-n)}{(1-n)[1-\delta+\delta a(\tau)]^2}\biggr)\biggr]^2\notag\\
&-\biggl[-\omega n[H-L]\Phi(\tau)\biggl(\frac{1-n+na(\tau)}{(1-n)[1-\delta+\delta a(\tau)]}\biggr)\notag\\
&+2\biggl(\widetilde{\omega}+\omega n[H-L]a(\tau)\biggr)\biggl(\frac{\Phi(\tau)(\delta-n)}{(1-n)[1-\delta+\delta a(\tau)]^2}\biggr)\notag\\
&+\biggl(\tau\widetilde{\omega}+\omega n[H-L]A(\tau)\biggr)\frac{(\delta-n)\Phi(\tau)[(\theta+1)(\delta-1)+\delta a(\tau)(\theta-1)]}{(1-n)(\tau+c)[1-\delta+\delta a(\tau)]^3}\biggr]\Biggr\}.
\end{align}
The previous formula has been obtained using the following results:
\begin{align*}
\frac{da(\tau)}{d\tau}&=\frac{d}{d\tau}\biggl[1-\int_0^{\tau}\Phi(x)dx\biggr]\\
&=-\Phi(\tau);
\end{align*}
\begin{align*}
\frac{d}{d\tau}\biggl[\frac{1-n+na(\tau)}{(1-n)[1-\delta+\delta a(\tau)]}\biggr]=&\frac{d}{d\tau}\biggl[\frac{\delta-n+n[1-\delta+\delta a(\tau)]}{\delta(1-n)[1-\delta+\delta a(\tau)]}\biggr]\\
=&\frac{\delta-n}{\delta(1-n)}\frac{d}{d\tau}\biggl[\frac{1}{1-\delta+\delta a(\tau)}\biggr]\\
=&\frac{\Phi(\tau)(\delta-n)}{(1-n)[1-\delta+\delta a(\tau)]^2};
\end{align*}
in the end
\begin{align*}
\frac{d^2}{d\tau^2}&\biggl[\frac{1-n+na(\tau)}{(1-n)[1-\delta+\delta a(\tau)]}\biggr]=\frac{\delta-n}{1-n}\frac{d}{d\tau}\bigg[\frac{\Phi(\tau)}{[1-\delta+\delta a(\tau)]^2}\biggr]\\
=&\frac{\delta-n}{(1-n)[1-\delta+\delta a(\tau)]^4}\Biggl\{-\Phi(\tau)\frac{\theta+1}{\tau+c}[1-\delta+\delta a(\tau)]^2\\
&+2\Phi(\tau)[1-\delta+\delta a(\tau)]\delta\Phi(\tau)\Biggr\}\\
=&\frac{(\delta-n)\Phi(\tau)}{(1-n)[1-\delta+\delta a(\tau)]^3}\Biggl\{\frac{\theta+1}{\tau+c}[-1+\delta-\delta a(\tau)]+2\delta\Phi(\tau)\Biggr\}\\
=&\frac{(\delta-n)\Phi(\tau)}{(1-n)(\tau+c)[1-\delta+\delta a(\tau)]^3}\biggl[(\theta+1)(\delta-1)-\delta\theta a(\tau)-\delta a(\tau)+2\delta\theta a(\tau)\biggr]\\
=&\frac{(\delta-n)\Phi(\tau)[(\theta+1)(\delta-1)+\delta a(\tau)(\theta-1)]}{(1-n)(\tau+c)[1-\delta+\delta a(\tau)]^3}.
\end{align*}
The last computation has been got since, recalling equation~\eqref{eqn:funa} and the one~\eqref{eqn:omori} in Appendix~\ref{sec:appc}, we have
\begin{align*}
\frac{d\Phi(\tau)}{d\tau}&=\theta c^{\theta}\frac{d}{d\tau}(\tau+c)^{-\theta-1}\\
&=-\theta c^{\theta}(\theta+1)(\tau+c)^{-\theta-2}\\
&=-\Phi(\tau)\frac{\theta+1}{\tau+c}
\end{align*}
and
\[
\Phi(\tau)=\frac{\theta a(\tau)}{\tau+c}.
\]

\section*{Notation}
\begin{tabular}{p{4cm} p{10.9cm}}
$a$, $\alpha$ & exponent of the productivity law (equation~\eqref{eqn:produttivita}, Appendix~\ref{sec:appc}); it holds $a=\alpha\ln10$.\\
$a(t)$ & auxiliary function defined in~\eqref{eqn:funa} as the integral of the Omori-Utsu law (equation~\eqref{eqn:omori}, Appendix~\ref{sec:appc}) between $t$ and infinity.\\
$A(t)$ & auxiliary function defined in~\eqref{eqn:funagrande} as the integral of $a(x)$ between zero and $t$.\\
$\beta$, $b$ & exponent of the Gutenberg-Richter law (equation~\eqref{eqn:GR}, Appendix~\ref{sec:appc}); $b$ is the so-called \emph{b-value} and holds that $\beta=b\ln10$.\\
$b(t)$ & auxiliary function defined in~\eqref{eqn:b} as the integral of the Omori-Utsu law (equation~\eqref{eqn:omori}, Appendix~\ref{sec:appc}) between zero and $t$. It holds $b(t)=1-a(t)$.\\
$c$ & parameter of the Omori-Utsu law (equation~\eqref{eqn:omori}, Appendix~\ref{sec:appc}).\\
$C_1$, $C_2$ & constants entering in the definition of the function $q(m')$; it holds $C_2=2C_1$.\\
$D(z;t,\tau,m_0,\overline{m})$ & auxiliary function defined in~\eqref{eqn:b19a} entering in the definition of $y(z;t,\tau,m_0,m')$.\\
$D_+(z;t,\widetilde{m},\overline{m})$ & auxiliary function defined in~\eqref{eqn:b19c} entering in the definition of $y(z;t,\tau,m_0,m')$ and $s(z;t,m_0,m')$.\\
$\delta$ & parameter defined in~\eqref{eqn:mm09}.\\
$\Delta$ & parameter defined in~\eqref{eqn:deltagrande}.\\
$f(m',m'')$ & auxiliary function defined in~\eqref{eqn:f3} (Appendix~\ref{sec:appa}) useful to calculate the conditional probability density $p(m''|m')$.\\
$f(z,m',m)$ & function defined in~\eqref{eqn:perz} useful to calculate both the probability generating function of the total number of events in $[0,\tau]$ and the one relative to the number of event triggered by a generic shock in the same time interval. The multiplication by this function allows us to add the triggering event if its magnitude is larger than the observability threshold $m$.\\
$\phi(t)$ & Omori-Utsu law (equation~\eqref{eqn:omori}, Appendix~\ref{sec:appc}); it can be interpreted as the probability generating function of random times at which first generation shocks independently occur, triggered by an event occurred in $t=0$.\\
$F_{\tau}(\tau)$ & inter-event time density.\\
$G(z;t,\tau,m')$ & probability generating function of the number of observable events triggered in $[0,\tau]$ by a spontaneous shock of magnitude $m'$ occurred in $-t$, with $t\ge0$.\\
$H$ & parameter defined in~\eqref{eqn:mm04}.\\
$I_A(m'')$ & magnitude integral from the reference cutoff $m_0$ to infinity, defined in~\eqref{eqn:Agrande}, of the product between the Gutenberg-Richter law~\eqref{eqn:GR} (Appendix~\ref{sec:appc}), the productivity law~\eqref{eqn:produttivita} (Appendix~\ref{sec:appc}) and the conditional probability density function~\eqref{eqn:pcond} (Appendix~\ref{sec:appa}) of triggered events' magnitude. We impose that $I_A(m'')=np(m'')$.
\end{tabular}

\begin{tabular}{p{4cm} p{10.9cm}}
$I_B(m'')$ & magnitude integral from the completeness value $m$ to infinity, defined in~\eqref{eqn:Bgrande}, of the product between the Gutenberg-Richter law~\eqref{eqn:GR} (Appendix~\ref{sec:appc}), the productivity law~\eqref{eqn:produttivita} (Appendix~\ref{sec:appc}) and the conditional probability density function~\eqref{eqn:pcond} (Appendix~\ref{sec:appa}) of triggered events' magnitude. We impose that $I_B(m'')\approx nLp(m'')$.\\
$\kappa$ & multiplicative parameter of the productivity law (equation~\eqref{eqn:produttivita}, Appendix~\ref{sec:appc}).\\
$L$ & parameter entering in the approximate condition that the integral $I_B(m'')$ must satisfy; it has to be found numerically.\\
$L(z;\tau)$ & auxiliary function useful to compute the probability generating function $\Omega(z;t)$ and then the probability to have zero events in $[0,\tau]$.\\
$\lambda$ & rate of the whole process of ''observable'' events.\\
$m$ & threshold value of completeness magnitude.\\
$m_0$ & threshold value of reference magnitude.\\
$n$ & average branching ratio defined alternatively in~\eqref{eqn:enne}.\\
$\overline{N}(t,\tau), \overline{N}_{m=m_0}(t,\tau)$, $\overline{N}_-(t)$ & auxiliary functions useful to compute $L(0;\tau)$; it holds $\overline{N}(t,\tau)=\overline{N}_{m=m_0}(t,\tau)$.\\
$\omega$ & rate corresponding to the ''observable'' events of the spontaneous component of the process.\\
$\widetilde{\omega}$ & parameter defined in~\eqref{eqn:omegam}.\\
$\Omega(z;\tau)$ & probability generating function of the total number $N(\tau)$ of observable events in $[0,\tau]$.\\
$p(m)$ & Gutenberg-Richter law (equation~\eqref{eqn:GR}, Appendix~\ref{sec:appc}); it is the probability density function of events' magnitude when we do not consider past seismicity.\\
$\mathbb{P}(\tau)$ & probability to have zero events in $[0,\tau]$.\\
$p(m''|m')$ & conditional probability density function assumed for triggered events' magnitude; it is defined in~\eqref{eqn:pcond} (Appendix~\ref{sec:appa}).\\
$p(m''|tr)$ & conditional probability for an event to have magnitude $m''$ given the fact that it is triggered; it is defined in~\eqref{eqn:pperfigli}.\\
PGF & Probability Generating Function.\\
$\Psi(h,\widetilde{m})$ & auxiliary function defined in~\eqref{eqn:psi} useful to compute $L(z;\tau)$.\\
$Q$ & parameter defined in~\eqref{eqn:Qgrande} and calculated as the integral of the Gutenberg-Richter law between the completeness value $m$ and infinity.\\
$q(m')$, $s(m')$ & auxiliary functions entering in the definition of $f(m',m'')$ (equation~\eqref{eqn:f3}); the function $q(m')$ is defined in~\eqref{eqn:q2} and holds that $s(m')=2q(m')$.\\
$\rho(m)$ & productivity law (equation~\eqref{eqn:produttivita}, Appendix~\ref{sec:appc}); it is the contribution of first generation shocks triggered by an event with magnitude $m'$.\\
$s(z;t,m_0,m')$, $y(z;t,\tau,m_0,m')$ & auxiliary functions used to define $L(z;\tau)$.\\
$\theta>0$ & exponent of the Omori-Utsu law (equation~\eqref{eqn:omori}, Appendix~\ref{sec:appc}).
\end{tabular}





\clearpage
\addcontentsline{toc}{section}{\refname}
\nocite{*}
\bibliographystyle{plain}
\bibliography{biblioteor}

\end{document}